\documentclass[reqno,10pt]{amsart}\usepackage{latexsym,amssymb,amsthm,amsmath,amscd,a4wide}

\newcommand{\M}{\ensuremath{\widetilde{M}}}

\theoremstyle{plain}
\newtheorem{theorem}{Theorem}[section]
\newtheorem{lemma}{Lemma}[section]
\newtheorem{proposition}{Proposition}[section]
\newtheorem{corollary}{Corollary}[section]

\theoremstyle{remark}
\newtheorem{remark}{Remark}[section]
\theoremstyle{remark}

\theoremstyle{definition}

\numberwithin{equation}{section}

\def\<{\left< }
\def\>{\right> }
\def\E4{\mathbb E^4 }

\def\E4{\mathbb E^4 }
\def\({$}
\def\){$}
\def\[{$$}
\def\]{$$}

\begin{document}


\title[Solitonic Inequalities in Riemannian Geometry]{A Comprehensive Review of Solitonic Inequalities in Riemannian Geometry}

\author[B.-Y. Chen, M. A. Choudhary, M. Nisar and M. D. Siddiqi]{Bang-Yen Chen, Majid Ali Choudhary, Mohammed Nisar\\ and Mohd Danish Siddiqi}




\begin{abstract}
In Riemannian geometry, Ricci soliton inequalities are an important field of study that provide profound insights into the geometric and analytic characteristics of Riemannian manifolds. An extensive study of Ricci soliton inequalities is given in this review article, which also summarizes their historical evolution, core ideas, important findings, and applications.
We investigate the complex interactions between curvature conditions and geometric inequalities as well as the several kinds of Ricci solitons, such as expanding, steady, and shrinking solitons. We also go over current developments, unresolved issues, and possible paths for further study in this fascinating area.
\end{abstract}

\maketitle

\section{Introduction}

The following equation was first studied in 1981 by R. Hamilton \cite{Ha82}:
\begin{equation}\label{Ricci1}
\frac{d}{d t} g(t)=-2 \operatorname{Ric}_{g(t)}, g_0 = g(0),
\end{equation}
where Ric is the Ricci tensor.

On a given manifold $M$, an unknown differentiable curve $t \mapsto g(t)$ of Riemannian metrics is subjected to a condition. It is commonly stated that the Ricci flow is produced by the solution to equation \eqref{Ricci1}. It is necessary for the derivative of the curve with respect to $t$ to equal twice the Ricci tensor of the metric $g(t)$ at each $t$ in its domain interval.

Specifically, define the curves $t \mapsto g(t)$ that are solutions (trajectories) of the Ricci flow, given an initial metric $g(0)$ and a maximum interval $[0, T)$ of the variable $t$, where $0\leqslant<t<\infty$.

In local coordinates $x^{j}$, the system of nonlinear PDE affects the component $g_{j k}=g\left(e_{j}, e_{k}\right)$ of the metrics $g=g(t)$. The functions depend on variables $t$ in addition to $x^j$. A relatively complex coordinate variant of condition \eqref{Ricci1} exists as well:
$$
\frac{\partial g_{j k}}{\partial t}=-2 R_{j k}, \;\; \text { where } \;\; R_{j k}=\partial_{p} \Gamma_{j k}^{p}-\partial_{j} \Gamma_{p k}^{p}+\Gamma_{q p}^{p} \Gamma_{j k}^{p}-\Gamma_{j q}^{p} \Gamma_{p k}^{q}.
$$
The Christoffel symbols for the metric $g=g(t)$ are $\Gamma_{j k}^{p}$, as defined by 
$$ 2 \Gamma_{j k}^{p}=g^{p q}\left(\partial_{j} g_{k q}+\partial_{k} g_{j q}-\partial_{q} g_{j k}\right),$$ 
where $\left[g^{j k}\right]$ represents the matrix inverse of $\left[g_{j k}\right]$ at each point of the coordinate domain, with $g^{j k}$ denoting the contravariant components of the metric. The repeated indices are added over in the following formulas for $R_{j k}$ and $2 \Gamma_{j k}^{p}$ because the Einstein summing convention was applied.\\

\subsection{Ricci flow in the proof of Poincar\'e conjecture}\label{S1.1}
In \cite{Ha82}, Hamilton established the uniqueness and existence of a maximal Ricci flow trajectory on every compact manifold, with any given initial metric $g(0)$. Additionally, he tried to use this finding to support the Poincaré conjecture in three dimensions \cite{Poincare1905}.

"The Hamilton program" was his suggested outline of such a proof, and it comprised a particular set of actions. In 2002, Grigori Perelman \cite{Pe02,Pe03,Pe03.2} completed the final steps.

Perelman also demonstrated the considerably general Thurston geometrization conjecture for three-dimensional manifolds at the same time.

Surgery plays a key role in Perelman's work, as they are required when the Ricci flow encounters a singularity in limited time $(T<\infty)$. After the procedure, in a topologically simpler scenario, the Ricci flow is once again employed.

\subsection{Ricci solitons - "fixed points" of the Ricci flow}\label{S1.2}
A Riemannian metric $g_0 = g(0)$ is a Ricci soliton on a manifold $M$.
It develops in an inessential manner under the Ricci flow, so up to multiplications by positive constants (also known as "rescalings") and diffeomorphisms, all of the stages $g(t)$ overlap with $g(0)$. Put another way, under the equivalence relationship that was just mentioned, this kind of metric denotes a fixed point in the quotient of the space of metrics on $M$ of the Ricci flow.

The aforementioned idea is evident in the case of compact manifolds since there exists only one maximal Ricci flow trajectory with a given initial metric. A Ricci soliton, in the absence of the compactness condition, is a Riemannian metric $g$ on a manifold $M$ with a solution to equation \eqref{Ricci1}, which reflects an inessential development of the Riemannian metric and meets the initial condition $g(0)=g$.

\subsection{Applications}\label{S1.3}
In the study of the Ricci flow and differential geometry, Ricci solitons are significant entities. These self-similar solutions to the Ricci flow equation are important in comprehending the creation of singularities and the long-term behavior of the Ricci flow.\\

\noindent\textbf{Understanding Singularities in Ricci Flow:}
Near singularities, the behavior of the Ricci flow is commonly modeled by Ricci solitons. They provide insight into the possible shapes and structures that singularities can take \cite{hamilton1988ricci}.\\

\noindent \textbf{Example:}
 The Cigar soliton is a two-dimensional steady Ricci soliton with the metric:
  \[
  ds^2 = \frac{dx^2 + dy^2}{1 + x^2 + y^2}
  \]
  It models a steady soliton that represents the behavior near certain types of singularities in two dimensions.\\

\noindent\textbf{Geometric Analysis and Topology:}
Ricci solitons are used in the classification of manifolds in terms of their curvature properties. They are involved in the research of metrics with non-negative Ricci curvature and Einstein manifolds \cite{cao2010complete}.\\

\noindent \textbf{Example:}
 Cylindrical Soliton: The product of a round sphere with Euclidean space, such as $( S^{n-1} \times \mathbb{R})$, is an example of a steady Ricci soliton. It models the behavior of Ricci flow near cylindrical singularities.\\

\noindent\textbf{Mathematical Physics and String Theory:}
Ricci solitons are used in string theory and the study of the renormalization group flow. They model fixed points and other critical phenomena in the theory of quantum fields \cite{streets2014symplectic}.\\
  
\noindent\textbf{Analysis of Heat Kernels and Eigenvalues:}
 Ricci solitons are used in the study of heat kernels and eigenvalues of the Laplace-Beltrami operator on manifolds. This has applications in spectral geometry and the study of heat diffusion on manifolds \cite{zhang2009completeness}.\\
  
\noindent\textbf{Applications In general relativity:} 
In general relativity, the Ricci curvature tensor plays a crucial role in describing the curvature of spacetime in terms of Einstein equation \cite{einstein1915feldgleichungen}. Energy conditions are theoretical constraints that relate the energy content of matter and fields to the curvature of spacetime. One such condition involves the Ricci curvature tensor.\\

 The energy condition related to the Ricci curvature tensor is called the "weak energy condition" (WEC) \cite{hawking2023large}. It states that for any timelike vector  $V^{\mu}$, the contraction $R_{\mu \nu} V^{\mu} V^{\nu} $ is non-negative, where  $R_{\mu \nu}$ represents the components of the Ricci curvature tensor.
Mathematically, the weak energy condition can be expressed as:
$$ R_{\mu \nu} V^{\mu} V^{\nu} \geq 0 $$
for all timelike vector fields $V^{\mu}$. This condition effectively says that any observer's measurement of the energy density must be non-negative.\\

\textbf{Connections to Perelman's Work:} Grigori Perelman's work \cite{73} on the geometrization conjecture and the Ricci flow with surgery is strongly related to the research of Ricci soliton inequalities. Comprehending and improving the inequalities linked to Ricci solitons adds to the larger program of comprehending three-manifold geometry and topology.

\section{Preliminaries}\label{S2}

\subsection{Discovery of solitons}\label{S2.1}
A soliton is a type of nonlinear wave that possesses two distinct properties:

\hskip.06in
(1) Without changing its shape or speed, a localized wave continues to propagate.

\hskip.06in
(2) Localized waves maintain their identities and are stable against mutual collisions.\\

The first is a condition of solitary waves that has been known since the 19\textsuperscript{th} century in hydrodynamics. The second indicates that the wave possesses a particle-like quality. In contemporary physics, a suffix-on, like phonon and photon, is employed to denote the particle property. Using the particle attribute of a solitary wave, Zabusky and Kruskal \cite{106} called it a "soliton." 

The story of how soliton was discovered is fascinating and astounding. The Scottish engineer and scientist John Scott-Russell (1844) made the first on record observation of the single wave in 1834.

Scott-Russell is credited as coining the term "solitary wave". Scientists like Stokes, Boussinesq, and Rayleigh were interested by this occurrence, but no one had an empirical explanation until 1898, when two Dutch physicists named Korteweg and de Vries published their now-famous equation \cite{105} (referred to as $KdV$ equation).

\begin{equation}\label{Pr1}
u_{t}+\alpha u u_{x}+\beta u_{x x x}=0., \quad \alpha\; \text{and}\; \beta\;\;  \text {are constants.} 
\end{equation}
In this case, $x$ is the coordinate traveling at the velocity of a linearized wave, and $u(x, t)$ represents the height above mean sea level. Equation \eqref{Pr1} (1-soliton) has a solution for solitary waves:
\begin{equation}\label{Pr2}
u(x, t)=\frac{3 v}{\alpha} \operatorname{sech}^{2}\left[\vartheta\right] 
\end{equation}
with $\vartheta = \frac{1}{2} \sqrt{\frac{v}{\beta}}(x-v t)$.

Pasta, Fermi and Ulam analyzed the equilibrium state approach in an a 1-dimensional non linear lattice in 1955 \cite{FPUT55}. An ergodic system, or an equal distribution of energy among all of the modes, was predicted as a consequence of the linear system's nonlinear interactions among its normal modes. Results from a numerical analysis contradicted this idea. The system eventually reverts to its original configuration, despite the energy being distributed unevenly among all the modes .

In 1965, Zabusky and Kruskal used a numerical solution of the $KdV$ equation as a model for a nonlinear lattice to discover the recurrent occurrences. They also discovered an unanticipated $KdV$ equation feature. Sharply peaks appear on a waveform that was initially smooth. These pulse-waves flow through one another after impacts and travel nearly independently at steady speeds. According to a thorough examination, each pulse is a $\operatorname{sech}^{2}$-type solitary wave, as shown by \eqref{Pr2}, and these solitary waves exhibit stable particle behavior. And thus, the discovery of the soliton.

\subsection{Establishment of soliton concept}\label{S2.2}

\begin{itemize}
\item[{1.}] How can one verify the soliton's properties analytically? 

\item[{2.}] Why do particles like solitons have stability?

\item[{3.}]  Is soliton a phenomenon specific to the equation \eqref{Pr1}? 
\end{itemize}

\noindent  In order to construct the concept of solitons, the answers to these questions were essential.

 After the independent and dependent variables are appropriately scaled, equation \eqref{Pr1}  takes on its shape.
\begin{equation}\label{Pr3}
u_{x x x}+u_{t}-6 u u_{x}=0 .
\end{equation}
A stationary waveform is produced by nonlinear and dispersion effects, which are represented by the $KdV$ equation. In spite of mutual interactions, the equation's preserved parameters guarantee the stability of isolated waves. Given the infinite degrees of freedom of the field variable, the dynamic properties of the system are bounded by infinite conservation laws, permitting an arbitrary number of solitons and infinite conserved values.

The inverse scattering approach is used to study the fundamental properties of solitons. In 1967, Gardner, Kruskal, Greene and Miura developed a linear problem known as the eigenvalue problem \cite{GGK67}. The potential function $u(x, t)$ is the solution to equation \eqref{Pr3}:
\begin{equation}\label{Pr4}
u(x, t) \psi-\psi_{x x}=\lambda \psi .
\end{equation}
It is  shown that the Eigenvalue of $\boldsymbol{\lambda}$ does not rely on time when $u$ evolves as per \eqref{Pr3}. In quantum physics, equation \eqref{Pr4} is simply the Schrödinger equation. The issue is deemed scattering when, for a given $u$, the discrete eigenvalues $\lambda_{n}=-\eta_{n}^{2}$, the transmission coefficient $1 / a(k)$, the coefficient of reflection $b(k) / a(k)$, and so forth are determined. On the other hand, an inverse scattering problem is the possibility of a specific collection of scattering data, like $a(k), b(k), \lambda_{n}$, etc. Gelfand-Levitan and Marchenko solved the latter problem for equation \eqref{Pr4}. Further information suggests that the eigen-function $\psi(x, t)$'s time development is
\begin{equation}\label{Pr5}
\psi_{t}=-4 \psi_{x x x}+3 u_{x} \psi+6 u \psi_{x} .
\end{equation}
 In summary, the ${KdV}$ equation's solution is provided by
\begin{equation}\label{Pr7}
u(x, t)=2 \frac{\partial}{\partial x} K(x, x ; t) ,
\end{equation}
where
\begin{align}\label{Pr8}
& K(x, y ; t)+\int_{-\infty}^{x} K(x, z ; t) F(z+y ; t) \mathrm{d} z+F(x+y ; t)=0 ,
\end{align}
\begin{align}\label{Pr9}
& F(x ; t)=\sum_{n=1}^{N} c_{n}^{2}(t) \mathrm{e}^{\eta_{n} x}+\frac{1}{2 \pi} \int_{-\infty}^{\infty} \frac{b(k, t)}{a(k, 0)} \mathrm{e}^{-i k x} \mathrm{~d} k .
\end{align}
The Gelfand-Levitan equation is the integral equation \eqref{Pr8}. It is important to keep in mind that all of the calculations in the example above are linear issues. This resolves the $KdV$ equation's initial-value problem. A soliton is associated with each limited state that has a discrete eigenvalue. In particular, it is simple to obtain the $\mathrm{N}$-soliton solution corresponding to $\mathrm{N}$ bound states by solving the Gelfand-Levitan equation where the coefficient for reflection $r(k, 0)=b(k, 0) / a(k, 0)$ is identically zero (the reflectionless potential). Through the exact formulation of the $\mathrm{N}$-soliton solution \cite{104}, we can show that the soliton is stable against mutual collisions and that the collisions are pairwise, meaning that they only result in position changes in the solitons. The above-discussed solution methodology is called the inverse scattering method.

It is possible to solve the $KdV$ equation's initial-value problem. At the time, though, it seemed like a coincidence. Five years later, Zakharov and Shabat \cite{ZS73} addressed the nonlinear Schrödinger (NLS) problem by extending on the inverse scattering technique.\\
\begin{equation}
i \psi_{t}+\psi_{x x}+2|\psi|^{2} \psi=0 .
\end{equation}
then afterwards Wadati resolved the equation for modified $KdV$ $(mKdV)$ \cite{110}.
\begin{equation}
u_{t}+6 u^{2} u_{x}+u_{x x X}=0 .
\end{equation}
Furthermore, we now have more than 100 soliton equations after solving the Sine-Gordon equation \cite{MA74}.

We may precisely define the soliton using the inverse scattering method. Transferring the field variables to the scattering data is known as a canonical transformation. The scattering data space contains definitions for the action-angle variables. Because of this, the soliton equation is a fully integrable system, and the soliton is a system fundamental mode \cite{ZF71}.

The greatest advance in the theory of solitons is represented by the discovery of the inverse scattering method. The inverse scattering method can be thought of as the Fourier transformation expanded into nonlinear problems if the scattering data space is viewed as an extension of the momentum space. The Fourier transform's development allowed for the solution of the diffusion problem in 1811. It was refined into a cohesive approach for solving nonlinear evolution problems after more than 150 years.

\subsection{Ricci soliton}\label{S2.3}
A Ricci soliton is a complete Riemannian manifold $(M, g)$ in differential geometry if and only if there is a smooth vector field $V$ such that
$$
\operatorname{Ric}(g)-\lambda g+\frac{1}{2} \mathcal{L}_V g= 0.
$$
Here, $\lambda \in \mathbb{R}$ is the constant. In this case, the Lie derivative is represented by $\mathcal{L}$, and the Ricci curvature tensor by $\mathrm{Ric}$. $(M, g)$ is a gradient Ricci soliton if $f: M \rightarrow \mathbb{R}$ exists and $V=\nabla f$. The soliton equation is then formed as:
$$
\operatorname{Ric}(g)+\nabla^2 f=\lambda g \text {. }
$$
Note that the Einstein equation is the solution of the previously mentioned equations when either $V=0$ or $f=0$. This is why  Einstein manifolds are the generalization of Ricci solitons.\\

\noindent \textbf{Role of \( \mathbf{V} \) in Ricci solitons} \cite{cao2010complete,chow2023hamilton}\\

\noindent \textbf{1. Characterization of soliton Type:}
\begin{itemize}
    \item Gradient Solitons: For any smooth function \( f \), \( \mathbf{V} = \nabla f \) is a gradient vector field. In this case, the Ricci soliton equation reduces to \( \text{Ric} + \nabla \nabla f-\lambda g =0  \). The potential function in this case is denoted by \( f \).
   
   \item Rotational Solitons: For rotational solitons, \( \mathbf{V} \) generates a 1-parameter family of isometries, implying certain symmetries in the manifold's structure.

   \item Homogeneous Solitons: In homogeneous solitons, \( \mathbf{V} \) exhibits a high degree of symmetry, often reflecting the underlying symmetry group of the manifold.
\end{itemize}

\noindent \textbf{2. Geometric Interpretation:}
   \( \mathbf{V} \) governs the self-similar evolution of the soliton under the Ricci flow. It determines how the metric \( g \) changes along the flow direction induced by \( \mathbf{V} \).\\

\noindent \textbf{3. Physical Interpretation:}
   - \( \mathbf{V} \) can be interpreted as a flow vector field that describes how the metric evolves to maintain the soliton condition. Depending on the sign of \( \lambda \), \( \mathbf{V} \) either contracts, expands, or keeps the metric unchanged (in the case of steady solitons).

\subsection{Self-similar Ricci flow solutions}\label{S2.4}
The Ricci flow equation has a self-similar solution that is produced by a Ricci soliton $\left(M, g_0\right)$.
$$
\partial_t g_t=-2 \operatorname{Ric}\left(g_t\right) .
$$
Specifically allowing
$$
\sigma(t):=1-2 \lambda t ,
$$
and integrating the field of vectors $X(t):=\frac{1}{\sigma(t)}$, which is time-dependent. Using $\Psi_0$ as the identity and $V$ to produce a collection of diffeomorphisms $\Psi_t$, one may obtain a Ricci flow solution $\left(M, g_t\right)$ by taking
$$
g_t=\sigma(t) \Psi_t^*\left(g_0\right) .
$$
The metric $g_0$ is pulled back by the diffeomorphism $\Psi_t$, as expressed in the following statement: $\Psi_t^*\left(g_0\right)$. Therefore, a Ricci soliton homothetically declines based on the sign of $\lambda$, stays constant, or expands up to diffeomorphism under Ricci flow.

\subsection{Types of Ricci soliton with example}\label{S2.5}
\textbf{(i) Shrinking $(\lambda>0) \quad$} 

\begin{itemize}
\item  Einstein manifolds with positive scalar curvature \cite{besse2007einstein}.
\item Gradient compact Kahler-Ricci shrinkers \cite{Cao96,Koiso90,WZ04}.
\item  The 4-dimensional BCCD shrinker \cite{BCCD24}.
\item  The 4-dimensional FIK shrinker \cite{FIK03}.
\item Round shrinking cylinder $S^{n-1} \times \mathbb{R}, n \geq 3$ \cite{hamilton1993formations}.
\item  Round shrinking  sphere $S^n, n \geq 2$ \cite{chow2023hamilton}.
\item  Shrinking Gaussian soliton $\left(\mathbb{R}^n, g_{\text {eucl }}, f(x)=\frac{\lambda}{2}|x|^2\right)$ \cite{hamilton1993formations}.
\end{itemize}

In a shrinking Ricci soliton, the metric contracts under the Ricci flow, resembling a shrinking of the manifold. A classic example is the Gaussian soliton on Euclidean space. Consider Euclidean space with the metric:
   \[
   g = dx^2 + dy^2 + dz^2.
   \]
   The Gaussian soliton is given by the metric:
   \[
   g = (1 + 4t)^{-\frac{1}{2}}(dx^2 + dy^2 + dz^2).
   \]
   where \( t \) is time. This soliton shrinks isotropically as time progresses.\\
   
\noindent \textbf{(ii) Steady $(\lambda=0)$} 
\begin{itemize}
\item The 2-dimensional cigar soliton
$\left(\mathbb{R}^2, g=\frac{d x^2+d y^2}{1+x^2+y^2}, V=-2\left(x \frac{\partial}{\partial x}+y \frac{\partial}{\partial y}\right)\right).$
\item The Bryant soliton with rotational symmetry of $3 d$ and its extension to higher dimensions \cite{Bry05}.
\item  Ricci flat manifold.
\end{itemize}
The  sphere with its standard metric is among the most basic examples. We will use the standard round metric derived from the ambient Euclidean space to represent the \(n\)-dimensional sphere immersed in \(\mathbb{R}^{n+1}\) as \(\mathbb{S}^n\).

The standard metric on \(\mathbb{S}^n\) is given by
$$ g = d\beta_1^2+ \cos^2\beta_1 d\beta_2^2+\cdots + \cos^2\beta_1 \cdots \cos^2\beta_{n-1}d\beta_n^2. $$
where $\{\beta_1,\ldots,\beta_{n}\}$ is a standard spherical coordinate system.

The metric \(g\) evolves in accordance with the following equation under the Ricci flow:
\[ \frac{\partial g}{\partial t}+2 \text{Ric}  =0,\]
where \(\text{Ric}\) is the Ricci curvature tensor. The metric itself determines the Ricci curvature of a round sphere, which is \(\text{Ric} = \frac{2(n-1)}{n} g\). As a result, the equation for Ricci flow becomes:
\[ \frac{\partial g}{\partial t} = -\frac{4(n-1)}{n} g ,\]
Now, consider a time-dependent metric \( g(t) \) on \(\mathbb{S}^n\) given by:
\[ g(t) = e^{4(n-1)t/n} g(0) . \]

The Ricci flow equation is satisfied by this metric, which may be directly confirmed. Furthermore, we may retrieve the initial metric \( g(0) \) if we set \( t = 0 \). As a result, \( g(0) \) is implied to be a stable Ricci soliton and a fixed point of the Ricci flow.\\

\noindent \textbf{(iii) Expanding $(\lambda<0) \quad$} 
\begin{itemize}
\item Negative scalar curvature Einstein manifolds: K\"ahler-Ricci solitons expansion on the complex line bundles $O(-k), k>n$ over $\mathbb{C} P^n, n \geq 1$.
\end{itemize}

In an expanding Ricci soliton, the metric expands under the Ricci flow. An example is the Bryant soliton\cite{morgan2007ricci}. Consider the Euclidean plane with the metric:
$$  g = dx^2 + e^{2x}(dy^2 + dz^2) .
$$
This metric represents a soliton with positive curvature and expansion in the \( x \)-direction.\\

\subsection{Different inequalities in Riemannian geometry}\label{S2.6}
Inequalities play a crucial role in the study of Riemannian manifolds, which are geometric spaces that generalize the concept of curved surfaces. Here are several key reasons why inequalities are important in this context:\\

\noindent \textbf{1. Geometric Inequalities:}
   Based on the Bonnet-Myers theorem \cite{jost2008riemannian}, we can say that \( M \) is compact if and only if its sectional curvature is bounded from below by a positive constant \( \kappa \). When represented mathematically, this is:
     \[ \text{if } \text{sec}(M) \geq \kappa \text{, then } M \text{ is compact.} \]
     Here, the sectional curvature  of \( M \) is denoted by \( \text{sec}(M) \).\\

\noindent 2. \textbf{Topological Inequalities:}
The Bishop-Gromov inequality provides a bound on a Riemannian manifold's diameter in terms of its sectional curvature \cite{cheeger1975comparison}. It states explicitly that for any two points $p$ and $q$ in a Riemannian manifold \( M \), we have
     \[ \text{dist}(p, q) \leq \frac{\pi}{\sqrt{\text{sec}(M)}},  \]
     where the distance between points \( p \) and \( q \) is indicated by \( \text{dist}(p, q) \).\\

\noindent 3. \textbf{Rigidity Theorems:}
   The sphere theorem is one instance of a rigidity theorem \cite{Ha82}. This means that, provided \( \kappa > 0 \) is true, a full, simply connected Riemannian manifold \( M \) that has a constant \( \kappa \) bounded the sectional curvature from above is isometric to the standard sphere \( S^n \) with constant curvature \( \kappa \).\\

\noindent 4. \textbf{Gradient Inequalities:}
    The Laplacian of a function$f$ and the gradient's norm on a Riemannian manifold are related by the Bochner inequality \cite{schoen1994lectures}:
     \[ \| \nabla^2 f \| \geq \frac{1}{n} \| \text{Hess}(f) \|^2, \]
     where \( \nabla^2 f \) denotes the Hessian of \( f \), and \( \text{Hess}(f) \) denotes the Hessian operator.\\

\noindent 5. \textbf{Isoperimetric Inequalities:}
A basic finding in geometry that establishes a relationship between the area and the perimeter of a closed surface or curve in Euclidean space is known as the isoperimetric inequality. 

Let $A$ be the area enclosed by a closed curve or surface in Euclidean space and let \( P \) be the perimeter or circumference of the same curve or surface. Then, the isoperimetric inequality states that for any such curve or surface,
$$4\pi A \leq P^2 .$$

In two dimension, equality is maintained only when the curve or surface is a circle, and in three dimensions, it is a sphere.
In other words, among all closed curves or surfaces with a given area, the one with the smallest perimeter or circumference is the circle (or sphere), and the ratio of its perimeter (or circumference) to its area is \( \frac{P}{A} = \frac{1}{r} \), where \( r \) is the radius of the circle (or sphere).

According to Cheeger's inequality, the first eigenvalue of the Laplacian operator has a lower constraint in terms of the isoperimetric constant of a Riemannian manifold.\\

\noindent 6. \textbf{Functional Inequalities:}
    A function \( f \) on a Riemannian manifold \( M \) has a constant \( C \) such that, according to the Poincaré inequality,
     \[ \| f - \bar{f} \|_2 \leq C \| \nabla f \|_2 , \]
     where \( \bar{f} \) represents the average of \( f \), the gradient of \( f \) is represented by \( \nabla f \) , and \( \| \cdot \|_2 \) represents the \( L^2 \) norm.

\subsection{Definitions}\label{S2.7}

\textbf{$D$-Homothetic Deformation:} 
\cite{blair2010,tanno1989,kon1985} 
Scaling a manifold's metric tensor by a positive constant is commonly referred to as a homothetic deformation. A homothetic deformation by a constant $c > 0$ yields a new metric $\tilde{g}$ defined as follows if $g$ is the metric tensor on a manifold $M$ such that
$ \tilde{g} = c^2 g .$ 
This idea is extended by D-homothetic deformation, which applies a differential scale that can change over the manifold. It is frequently examined in relation to specific geometric structures, including paracontact or virtually contact metric structures.

They explored the behavior of virtually contact metric structures which are common in contact geometry under scaling transformations using D-homothetic deformations. Assume that a manifold $M$ has a metric structure $(\phi, \xi, \eta, g)$ which is almost contact . This case has the following: a vector field (the Reeb vector field), a 1-form (named $\eta$, a metric tensor (called $g$, and a (1,1)-tensor field called $\phi$.

A $D$-homothetic deformation involves transforming these components :
\[ \quad \tilde{\xi} = \frac{1}{c} \xi, \tilde{\phi} = \phi, \quad \tilde{\eta} = c \eta, \quad \tilde{g} = gc^2  + (c^2 - 1) \eta \otimes \eta, \]
where \( c \) is a positive constant.\\

\noindent\textbf{Conformal Ricci flow equations:}
\cite{fischer2004introduction} A modification to the conventional Ricci flow equation that yields, rather than the unit volume restriction, a scalar curvature constraint. Because of the part that conformal geometry plays in limiting scalar curvature, the resulting equations are also known as conformal Ricci flow equations since they are the vector field sum of a Ricci flow equation and a conformal flow equation. The formulas are given by
$$
\begin{aligned}
2\left(\frac{1}{n}g+\operatorname{Ric} g\right)+\frac{\partial g}{\partial t} & =-p g,\quad
R(g)  =-1,
\end{aligned}
$$
 for a metric $g$ that is dynamically changing and a static scalar field $p$. The conformal Ricci flow equations and the Navier-Stokes equations are the same in fluid mechanics.
$$
\begin{aligned}
\nu \Delta v+\nabla_{v} v+\operatorname{grad} p+\frac{\partial v}{\partial t} & =0,\quad 
\operatorname{div} v  =0.
\end{aligned}
$$
\noindent This comparison results in the time-dependent scalar field $p$ and the term "conformal pressure". Like the real physical pressure in fluid mechanics, which maintains the fluid's incompressibility, the conformal pressure acts as a Lagrange multiplier to conformally adjust the metric flow in order to preserve the scalar curvature constraint.\\

\noindent Conformal Ricci flow equations with the Einstein constant $-\frac{1}{n}$ have Einstein metrics at their locations of equilibrium. The formula $-2\left(\frac{1}{n}(g)+\operatorname{Ric} g\right)$ thus functions as a nonlinear restorative force and quantifies the flow's divergence from an equilibrium point. Unless it is in an equilibrium location, the conformal pressure $p \geq 0$ is positive. The constraint force $-p g$ conformally deforms $g$ to maintain the scalar curvature by acting pointwise orthogonally to the nonlinear restoring force $-2\left(\frac{1}{n}(g)+\operatorname{Ric} g\right)$.\\

\noindent \textbf{Conformal Ricci soliton:} \cite{chow2023hamilton,cao2009recent} A conformal Ricci soliton generalizes Ricci soliton concept by allowing a conformal change in the metric. Formally, a conformal Ricci soliton is defined by
\[ \operatorname{Ric}(g) + \mathcal{L}_X g = \lambda g + \mu g \]
where \(\mu\) is a differentiable function on \(M\).\\

\noindent This equation indicates that the Ricci curvature, modified by the Lie derivative of the metric, is equal to a scalar multiple of the metric plus a term that varies smoothly over the manifold.\\

\noindent \textbf{Conformal $\eta$ Ricci soliton:}
In 2015, A. Bhattacharyya and N. Basu \cite{basu2015conformal} introduced the conformal type Ricci soliton, which is a generalization of the classic type Ricci soliton and is represented by the equation
\begin{equation}\label{eta1}
\mathcal{L}_{V} g+2 \text { Ric }=\left[2 \lambda-p-\frac{2}{n}\right] g, 
\end{equation}
where $\lambda$ represents a constant and $p$ the conformal pressure. It is crucial to keep in mind that the conformal type Ricci soliton contains a self-similar solution to Fisher's conformal type Ricci flow equation. Following this, a multitude of investigations into conformal type Ricci solitons have been carried out on a range of geometric structures, including $(L C S)_{n}$-manifolds \cite{ganguly2020study} and generalized Sasakian space forms \cite{ganguly2021conformal}.

Gradient Ricci soliton refers, once more, to the situation in which the potential vector field $V$ in \eqref{eta1} is the gradient of a smooth function $f$ on $M$. The Ricci soliton potential function is the name given to this function, $f$. J. T. Cho and M. Kimura \cite{cho2009ricci} introduced the $\eta$-Ricci soliton in complex space forms for the first time. Later, C. Calin and M. Crasmareanu \cite{calin2012eta} studied it on Hopf hypersufaces. A $\eta$-Ricci soliton is said to be accepted by the Riemannian manifold $(M, g)$ when, for a smooth vector field $V$, the metric $g$ satisfies the following equation:
\begin{equation}\label{eta2}
2 \mu \eta \otimes \eta+\mathcal{L}_{V} g+\mathcal{L}_{V} g+2 \lambda g+2S=0, 
\end{equation}
where $\mathcal{L}_{V}$ is the Lie derivative in the direction of $V, S$, the Ricci tensor, and $\lambda$, $\mu$, are real constants. Observe that when $\mu=0$, the $\eta$-Ricci soliton becomes a Ricci soliton.

The combination of equations \eqref{eta1} and \eqref{eta2}: The following equation represents the conformal $\eta$-Ricci soliton, which was established by M. D. Siddiqi \cite{siddiqi2018conformal}.
\begin{equation*}
2 \mu \eta \otimes \eta+\mathcal{L}_{V} g+\left[2 \lambda-p-\frac{2}{n}\right] g+2 S=0 ,
\end{equation*}
where $n$ is the manifold's dimension, $p$ is  scalar field which is non-dynamical, $\mathcal{L}_{V}$ is the Lie derivative along $V$ and  $S$ the Ricci tensor's direction.\\

\textbf{Hitchin--Thorpe inequality:} \cite{hitchin1974compact,thorpe1969some,donaldson1997geometry} This is a significant finding in differential geometry, particularly for the investigation of four-dimensional manifolds. It provides a link between the Euler characteristic and the signature of an oriented compact 4-manifold accepting a Riemannian metric with non-negative scalar curvature.

For a differentiable compact orientated 4-manifold, the signature \( \tau(M) \)  and the Euler characteristic \( \chi(M) \) of \( M \) are related  by the Hitchin--Thorpe inequality as follows:
\[ 3|\tau(M)|\leq2\chi(M) . \]
\textbf{Hyperbolic Ricci soliton:} \cite{azami2023hyperbolic}
\noindent If a vector field $X$ on $M$ has real scalars $\lambda$ and $\mu$, such that a Riemannian manifold $\left(M^{n}, g\right)$, then the structure is known as hyperbolic Ricci soliton structure:
\begin{equation*}\label{hrcn}
\operatorname{Ric}+\frac{1}{2} \mathcal{L}_{X}\left(\mathcal{L}_{X} g\right)+\lambda \mathcal{L}_{X} g=\mu g .
\end{equation*}

\noindent As previously stated, depending on whether $\lambda$ is positive, zero, or negative, a hyperbolic Ricci soliton is said to be expanding, stable, or shrinking. Furthermore, the parameter $\mu$ represents the rate of change of our solutions and has geometric meaning. The equation above reduces to Einstein equation when $X=0$ or $X$ is a Killing vector field.

 The gradient hyperbolic structure ($M, g, X, \lambda, \mu$) is a hyperbolic Ricci soliton structure if there exists a function $f$, which is often referred to as the potential function, such that $X=\nabla f$. Hence, \eqref{hrcn} becomes
$$
\text { Ric }+2 \lambda \operatorname{Hess} f+\mathcal{L}_{\nabla f} \text { Hess } f=\mu g.
$$

\noindent \textbf{Example:} Let $\left(M^{n}, g(t)\right)=\left(\mathbb{R}^{n}, g_{\mathrm{can}}\right)$ be a static Euclidean space, then it is a stationary solution to the hyperbolic Ricci flow as given below
$$
\frac{\partial^{2} g_{i j}}{\partial t^{2}}=-2 R_{i j}, \quad g(0)=g_{0}, \quad \frac{\partial g_{i j}}{\partial t}=0
$$
\noindent and is therefore a constant hyperbolic soliton. A hyperbolic Ricci soliton that expands or shrinks modulo diffeomorphisms could also be applied to this solution.\\

\section{Solitonic inequalities on Sasakian manifold}\label{S3}

We have collected the followings on solitonic Inequalities.

\subsection{Inequalities involving Ricci soliton  on para Sasakian manifold}\label{S3.1}

In 2022, Choudhary et al. \cite{CKD22} examined the scalar curvature for submanifolds of Ricci solitons in order to investigate the connection between extrinsic and intrinsic invariants. The Ricci tensor was considered by  
\begin{equation}
\tilde{S}\left(\ell_{1}, \ell_{2}\right)=\left.\tilde{S}\right|_{T_{p} M}\left(\ell_{1}, \ell_{2}\right)+\left.\tilde{S}\right|_{T_{p}^{\perp} M}\left(\ell_{1}, \ell_{2}\right)
\end{equation}
for any $  \ell_{1}, \ell_{2} \in T_{p} M$.

Now using the above result, they deduce the following.

\begin{lemma} \cite{CKD22}\label{L:3.1}
Let $(\widetilde{M}, g, V, \Lambda)$ be a Ricci soliton and $M$ a contact pseudo-slant submanifolds of a ($\epsilon$)- para Sasakian form $\widetilde{M}(k)$. We obtain
\begin{equation}\begin{aligned}
\operatorname{div}(V)-\mathcal{A}_{1}\,&+\|\sigma\|^{2}\left\{\epsilon(k+1)\left(-\epsilon \cos ^{2} \theta-2\right)+(n-1)(k-3)\right\}(n) \\&
+\mathcal{K}_{2}\left(\epsilon\left(2 \epsilon-\cos ^{2} \theta-1\right)+5-n\right) 
+n \Lambda=0 ,
\end{aligned}
\end{equation}
where $\mathcal{A}_{1}=n^{2}\|H\|^{2}$ and $ \mathcal{K}_{2}=\sum_{\alpha=n+1}^{m} \sum_{i, j \neq 1}^{n}\big(\sigma_{i j}^{\alpha}\big)^{2}
$.
\end{lemma}

\begin{theorem} \cite{CKD22}   
 Let $(\tilde{M}, g, V, \lambda)$ is a Ricci soliton, while $M$ is a contact pseudo-slant submanifold of an $(\epsilon)$-para Sasakian form $\widetilde{M}(k)$. We get.
\begin{align*}
\operatorname{div}(V) \leq \,& n(n-1)\|H\|^{2}-n \Lambda 
 -\left\{\epsilon(k+1)\left(-2-\epsilon \cos ^{2} \theta\right)+(n)(n-1)(k-3)\right\} \\
& -\mathcal{K}_{2}\left(5-n+\epsilon\left(2 \epsilon-1-\cos ^{2} \theta\right)\right) .
\end{align*}
\end{theorem}

 Blaga and Carasmareanu \cite{BC20} proposed an inequality for a lower boundary of the geometry of $g$ in terms of gradient Ricci soliton for a  function $\psi$ on ambient space $M$ such that
\begin{equation}\label{3.15}
\|H e s s\|_{g}^{2}-\frac{1}{n}(\Delta \psi)^{2}\leq \|S\|_{g}^{2}  ,
\end{equation}
where {\it Hess} denotes the Hessian of the differentiable function $\psi$ on $M$. Let the soliton vector field $V$ satisfy $V=\nabla \psi$. The authors used \eqref{3.15} to articulate the following theorem.

\begin{theorem}\cite{CKD22} \label{T:3.2}
Consider a gradient Ricci soliton $(\tilde{M}, g, \nabla \psi, \Lambda)$ with a gradient-type soliton vector field $V$. Let $M$ represent a contact pseudo-slant submanifold of a $\tilde{M}(k)$-para Sasakian form ($\epsilon$). We can write,
\begin{align*}
\|S\|_{g}^{2} \geq &\, \| \text {Hess}\left\|_{g}^{2}-n^{3}\right\| H\left\|^{4}+\right\| \sigma \|^{4}+n \Lambda^{2} 
 +\left\{\epsilon(k+1)\left(-2-\epsilon \cos ^{2} \theta\right)+(n-1)(k-3)\right\}^{2} \\& +\frac{1}{n}\left(\mathcal{K}_{2}\right)^{2}\left(5-n+\epsilon\left(2 \epsilon-1-\cos ^{2} \theta\right)\right)^{2} .
\end{align*}
\end{theorem}

\begin{theorem}\cite{CKD22}\label{T:3.3}
Let $M$ be a submanifold of a ($\epsilon$)-para Sasakian form $\tilde{M}(k)$, with complete umbilical contact pseudo-slant. Assume that \eqref{T:3.2} holds true in its entirety. Next,
\begin{align*}
\|S\|_{g}^{2} \geq &\, \| \text {Hess} \|_{g}^{2}+n \Lambda^{2} 
 +\left\{\epsilon(k+1)\left(-2-\epsilon \cos ^{2} \theta\right) \epsilon+(n-1)(k-3)\right\}^{2} \\
& +\frac{1}{n}\left(\mathcal{K}_{2}\right)^{2}\left(5-n+\epsilon\left(2 \epsilon-1-\cos ^{2} \theta\right)\right)^{2} . 
\end{align*}
\end{theorem}
In this article, the authors \cite{CKD22} discovered an unknown inequality for pseudo-slant submanifolds of the $ \epsilon $ -para Sasakian manifold, notably concerning gradient Ricci soliton and gradient type vector field.

\subsection{Sasakian manifold using generalized Wintgen inequalities}\label{S3.2}

In \cite{choudhary2024generalized}, The authors create an inequality for Ricci solitons in order to find relationships between extrinsic and intrinsic invariants. 
 They proved the following results involving inequalities of Ricci soliton.

\begin{theorem}\label{T:3.4}
Let $M$ be a submanifold of a $(\epsilon)-P S S F \widetilde{M}(k)$ adopting a Ricci soliton with a potential vector field $V \in T M$ of Ricci soliton. Assume that $(\tilde{M}, g, V, \Lambda)$ is a Ricci soliton. The Ricci soliton is therefore expanding, steady, and shrinking in accordance with

\vskip.06in
\noindent  1. $\frac{n}{2}\|H\|^{2}<\frac{1}{n}\left(\operatorname{div}(V)+\|\sigma\|^{2}+\tau\right) \text {, }$ 
\vskip.06in

\noindent 2. $\frac{n}{2}\|H\|^{2}=\frac{1}{n}\left(\operatorname{div}(V)+\|\sigma\|^{2}+\tau\right)$, and
\vskip.06in

\noindent 3. $\frac{n}{2}\|H\|^{2}>\frac{1}{n}\left(\operatorname{div}(V)+\|\sigma\|^{2}+\tau\right)$, respectively.
\end{theorem}

\begin{corollary}\label{C:3.1}
Let $M$ be a submanifold of a $(\epsilon)$-PSSF $\widetilde{M}(k)$, and let $(\tilde{M}, g, V, \Lambda)$ be a gradient Ricci soliton with a potential vector field $V=\nabla(\psi) \in T M$ of gradient type. Gradient Ricci soliton is then admitted by the submanifold $M$ of $(\epsilon)-P S S F \widetilde{M}(k)$ if and only if
\begin{equation*}
\frac{n}{2}\|H\|^{2}<\frac{1}{n}\left(\Delta(\psi)+\|\sigma\|^{2}+\tau\right) .
\end{equation*}
\end{corollary}

\begin{theorem}\label{T:3.5}
Let $M$ be a submanifold of a $(\epsilon)-P S S F$ $\widetilde{M}(k)$, and let $(\tilde{M}, g, V=\nabla(\psi), \Lambda)$ be a gradient shrinking Ricci soliton with $V=\nabla(\psi) \in T M$ of gradient type. Then
\begin{equation*}
|W| \leq \frac{2}{n\|H\|^{2}-2\|\sigma\|^{2}-2 n \Lambda-2 \Delta(\psi)} \sqrt{\frac{2(n-1)}{n-2}} .
\end{equation*}
\end{theorem}

\begin{theorem}\label{T:3.6}
Let $(\tilde{M}, g, V=\nabla(\psi), \Lambda)$ be a gradient Ricci soliton of gradient type, and let $M$ be a CPS-submanifold of a ( $\epsilon)$-PSSF $\widetilde{M}(k)$. Then
\begin{align*}
& |W| \leq\left\{\left(\frac{k-3}{4}\right)(2 p+q)+\left(\frac{k+1}{4}\right)\left(-2 \epsilon-\epsilon^{2} \cos ^{2} \theta\right)\right\}(2 p+q+1)  \\
& \quad +\left(\frac{k+1}{4}\right)\left(4-2 p-q+2 \epsilon^{2}-\epsilon-\epsilon \cos ^{2} \theta\right)+(2 p+q+1)^{2}\|H\|^{2}-\|\sigma\|^{2}.
\end{align*}
\end{theorem} 

\begin{theorem} \label{T:3.7}
Given $V=\nabla(\psi) \in T M$ of gradient type and $M$ a completely umbilical CPS-submanifold of a ( $\epsilon)-P S S F \widetilde{M}(c)$, let $(\tilde{M}, g, V=\nabla(\psi), \Lambda)$ be a gradient shrinking Ricci soliton. Then
\begin{align*}
& |W| \leq\left\{(2 p+q)\left(\frac{k-3}{4}\right)+\left(-2 \epsilon-\epsilon^{2} \cos ^{2} \theta \right)\left(\frac{k+1}{4}\right)\right\}(2 p+q+1)  \tag{5.12}\\& \quad
 +\left(4-2 p-q+2 \epsilon^{2}-\epsilon-\epsilon \cos ^{2} \theta\right)\left(\frac{k+1}{4}\right)+\|H\|^{2}(2 p+q+1)^{2}-\|\sigma\|^{2}.
\end{align*}
\end{theorem}

\section{Solitonic inequalities on Kenmotsu manifold}\label{S4}

\subsection{In a Da-homothetically deformed Kenmotsu manifold}\label{S4.1}
In 2022, a Da-homothetically deformed Kenmotsu manifold was investigated by A. M. Blaga \cite{Blaga22} for almost Riemann and almost Ricci solitons. A gradient vector field, a solenoidal vector field, or the Reeb vector field of the deformed structure are possible vector fields for these solitons. In some situations, it is possible to obtain the scalar and Ricci curvatures explicitly.

In \cite{Blaga22}, the author derived the following results involving Ricci soliton inequalities.

\begin{theorem}\cite{Blaga22}\label{T:4.1}
   In a $(2 n+1)$-dimensional Da-homothetically deformed Kenmotsu manifold $(M, \bar{\phi}, \bar{\xi}, \bar{\eta}, \bar{g})$,If a gradient almost Riemann soliton is defined by $(V=\overline{\operatorname{grad}}(f), \bar{\lambda})$, then
\begin{equation*}
|\overline{\operatorname{Ric}}|_{\bar{g}}^{2} \geq(2 n-1)^{2}\left[\upsilon\right],
\end{equation*}
where $\upsilon=|\overline{\operatorname{Hess}}(f)|_{\bar{g}}^{2}-\frac{1}{2 n+1}(\bar{\Delta}(f))^{2}$.
\end{theorem}

\begin{corollary}\cite{Blaga22}\label{C:4.1}
    In a $(2 n+1)$-dimensional Kenmotsu manifold $(M, \bar{\phi}, \bar{\xi}, \bar{\eta}, \bar{g})$ that is Da-homothetically deformed, if a solenoidal gradient practically Riemann soliton is defined by $(V=\overline{\operatorname{grad}}(f), \bar{\lambda})$, then
\begin{equation*}
|\overline{\operatorname{Ric}}|_{g}^{2} \geq(2 n-1)^{2}|\overline{\operatorname{Hess}}(f)|_{\bar{g}}^{2} .
\end{equation*}
\end{corollary}

\begin{proposition}\label{P:4.1}\cite{Blaga22}
    In a $(2 n+1)$-dimensional D-homothetically deformed Kenmotsu manifold $(M, \bar{\phi}, \bar{\xi}, \bar{\eta}, \bar{g})$, if $(V=\overline{\operatorname{grad}}(f), \bar{\lambda})$ defines a gradient almost Riemann soliton, then
$$(2 n-1)^{2}|\operatorname{Hess}(f)|_{g}^{2}-4 n \frac{a-1}{a} \mathrm{scal}-4 n^{2}(2 n+1)\left(\frac{a-1}{a}\right)^{2}$$
$$\begin{gathered}
-\frac{(2 n-1)^{2}}{2 n+1}(\Delta(f))^{2}-\frac{2(2 n-1)^{2}}{2 n+1} \frac{a-1}{a}[\xi(f)-\xi(\xi(f))] \Delta(f)\\
+\frac{2 n(2 n-1)^{2}}{2 n+1}\left(\frac{a-1}{a}\right)^{2}(\xi(f))^{2}-\frac{2(2 n-1)^{2}(n+n a+a)(a-1)}{(2 n+1) a^{2}}(\xi(\xi(f)))^{2} \\
+\frac{2(2 n-1)^{2}(2 n+a)(a-1)}{(2 n+1) a^{2}} \xi(f) \cdot \xi(\xi(f)) \leq|\operatorname{Ric}|_{g}^{2}.
\end{gathered}$$
\end{proposition} 

\begin{corollary}\cite{Blaga22}\label{C:4.1.2}
As per Proposition \ref{P:4.1}'s hypothesis, if $V$ is $\bar{g}$-orthogonal to $\xi$, then 
\begin{gather*}
(2 n-1)^{2}|\operatorname{Hess}(f)|_{g}^{2}-\frac{(2 n-1)^{2}}{2 n+1}(\Delta(f))^{2}+\frac{4 n(2 n-1)(a-1)}{a} \Delta(f)  \\
+\frac{4 n^{2}(2 n+1)\left(a^{2}-1\right)}{a^{2}} \leq |\operatorname{Ric}|_{g}^{2} .
\end{gather*}
Furthermore, if $f$ is a function that is $\Delta$-harmonic, then
$$
|\operatorname{Ric}|_{g}^{2} \geq(2 n-1)^{2}|\operatorname{Hess}(f)|_{g}^{2}+\frac{4 n^{2}(2 n+1)\left(a^{2}-1\right)}{a^{2}} .
$$
\end{corollary}

\begin{proposition}\label{P:4.2}\cite{Blaga22}
$(M, \bar{\phi}, \bar{\xi}, \bar{\eta}, \bar{g})$ is a $(2 n+1)$-dimensional $Da$-homothetically deformed Kenmotsu manifold.if a solenoidal gradient practically Riemann soliton is defined by $(V=\overline{\operatorname{grad}}(f), \bar{\lambda})$, then
$$
|\operatorname{Ric}|_{g}^{2} \geq(2 n-1)^{2}|\operatorname{Hess}(f)|_{g}^{2}+\frac{a^{2}-1}{a^{2}}\left[4 n^{2}-(2 n-1)^{2}(\xi(\xi(f)))^{2}\right] .
$$
\end{proposition}

\begin{theorem}\cite{Blaga22}\label{T:4.1.2}
In a $(2 n+1)$-dimensional Kenmotsu manifold $(M, \bar{\phi}, \bar{\xi}, \bar{\eta}, \bar{g})$ that is Da-homothetically deformed, if the gradient defined by $(V=\overline{\operatorname{grad}}(f), \bar{\lambda})$ is almost Ricci soliton, then
\begin{equation*}
|\overline{\operatorname{Ric}}|_{\bar{g}}^{2} \geq|\overline{\operatorname{Hess}}(f)|_{\bar{g}}^{2}-\frac{1}{2 n+1}(\bar{\Delta}(f))^{2} .
\end{equation*}
\end{theorem}

\begin{corollary}\cite{Blaga22}\label{C:4.1.3}
If $(V=\overline{\operatorname{grad}}(f), \bar{\lambda})$ defines a solenoidal gradient almost Ricci soliton in a (2n+1)-dimensional Da-homothetically deformed Kenmotsu manifold $(M, \bar{\phi},\bar{\xi}, \bar{\eta}, \bar{g})$, then
\begin{equation*}\label{90}
|\overline{\operatorname{Ric}}|_{\bar{g}}^{2} \geq|\overline{\operatorname{Hess}}(f)|_{\bar{g}}^{2} .
\end{equation*}\end{corollary}

\begin{proposition}\cite{Blaga22}\label{P:4.1.3}
In a $(2 n+1)$-dimensional Da-homothetically deformed Kenmotsu manifold $(M, \bar{\phi}, \bar{\xi}, \bar{\eta}, \bar{g})$, if $(V=\overline{\operatorname{grad}}(f), \bar{\lambda})$ defines a gradient almost Ricci soliton, then
\begin{gather*}
|\operatorname{Ric}|_{g}^{2} \geq|\operatorname{Hess}(f)|_{g}^{2}-4 n \frac{a-1}{a} \text { scal }-4 n^{2}(2 n+1)\left(\frac{a-1}{a}\right)^{2}-\frac{1}{2 n+1}(\Delta(f))^{2}  \\
-\frac{2}{2 n+1} \frac{a-1}{a}[\xi(f)-\xi(\xi(f))] \Delta(f)+\frac{(a-1)2(2 n+a)}{(2 n+1) a^{2}} \xi(f) \cdot \xi(\xi(f)) \\
\quad+\frac{2 n}{2 n+1}\left(\frac{a-1}{a}\right)^{2}(\xi(f))^{2}-\frac{(a-1)2(n+n a+a)}{(2 n+1) a^{2}}(\xi(\xi(f)))^{2} .
\end{gather*}
\end{proposition}

\begin{corollary}\cite{Blaga22}\label{C:4.1.4}
 As per Proposition \ref{P:4.2}'s hypotheses, if $V$ is $\bar{g}$-orthogonal to $\xi$, then we obtain
\begin{equation*}
|\operatorname{Ric}|_{g}^{2} \geq|\operatorname{Hess}(f)|_{g}^{2}-\frac{1}{2 n+1}(\Delta(f))^{2}+\frac{(a-1)4 n}{a} \Delta(f)+\frac{(2 n+1)\left(a^{2}-1\right)4 n^{2}}{a^{2}} .
\end{equation*}
Furthermore, if $f$ is a function that is $\Delta$-harmonic, then
$$|\operatorname{Ric}|_{g}^{2} \geq|\operatorname{Hess}(f)|_{g}^{2}+\frac{4 n^{2}(2 n+1)\left(a^{2}-1\right)}{a^{2}}.
$$ \end{corollary}

\begin{proposition}\cite{Blaga22}\label{P:4.1.4}
In a (2n+1)-dimensional Da-homothetically deformed Kenmotsu manifold  $(M, \bar{\phi}, \bar{\xi}, \bar{\eta}, \bar{g})$, if $(V=\overline{\operatorname{grad}}(f), \bar{\lambda})$ defines a solenoidal gradient almost Ricci soliton, then
$$
|\operatorname{Ric}|_{g}^{2} \geq|\operatorname{Hess}(f)|_{g}^{2}+\frac{a^{2}-1}{a^{2}}\left[4 n^{2}-(\xi(\xi(f)))^{2}\right] .
$$
\end{proposition}

\subsection{Indefinite Kenmotsu manifolds as a framework for conformal $ \eta $-Ricci solitons}\label{S4.2}
Y. L. Li et al. \cite{LGDB22} investigated certain unique varieties of Ricci tensor in relation to $\epsilon$-Kenmotsu manifolds' conformal $\eta$-Ricci soliton in 2022. Subsequently, the writer examined certain curvature circumstances that permit conformal $\eta$-Ricci solitons on $\epsilon$-Kenmotsu manifolds.

\begin{corollary}\cite{LGDB22}\label{C:4.2.1}
An $n$-dimensional $\epsilon$-Kenmotsu manifold $(M, g)$ admits a conformal Ricci soliton $(g, \xi, \lambda)$. Consequently, $(M, g)$ becomes a $\eta$-Einstein manifold, and $\lambda$ satisfies $\lambda=\left(\frac{p}{2}+\frac{1}{n}\right)+\epsilon(n-1)$. Moreover, we have:
\vskip.1in

1. In the event that $\xi$ is spacelike, the soliton is either expanding, stable, or contracting as, $\left(\frac{p}{2}+\frac{1}{n}\right)>(1-n)$, $\left(\frac{p}{2}+\frac{1}{n}\right)=(1-n)$ or $\left(\frac{p}{2}+\frac{1}{n}\right)<(1-n) ;$ and

\vskip.1in
2. In the event that $\xi$ is timelike, the soliton is expanding, steady, or shrinking as, $\left(\frac{p}{2}+\frac{1}{n}\right)>(n-1)$, $\left(\frac{p}{2}+\frac{1}{n}\right)=(n-1)$ or $\left(\frac{p}{2}+\frac{1}{n}\right)<(n-1)$.
\end{corollary}

\begin{corollary}\cite{LGDB22}\label{C:4.2.2}
$\lambda=\left(\frac{p}{2}+\frac{1}{n}\right)+\epsilon(n-1)$ describes the relationship between the scalars $\lambda$ and $\mu$ on the manifold $(M, g)$, which is a $\eta$-Einstein manifold. An $n$-dimensional $\epsilon$-Kenmotsu manifold $(M, g)$ is constant multiple of $\xi$ if it admits a conformal Ricci soliton $(g, V, \lambda, \mu)$ such that $V$ is pointwise collinear with $\xi$. Moreover, we have:
\vskip.1in

1. In the event that $\xi$ is spacelike, the soliton is either expanding, steady, or shrinking as, $\left(\frac{p}{2}+\frac{1}{n}\right)+n>1$, $\left(\frac{p}{2}+\frac{1}{n}\right)+n=1$ or $\left(\frac{p}{2}+\frac{1}{n}\right)+n<1 ;$ and
\vskip.1in

2. In the event that $\xi$ is timelike, the soliton is expanding, steady, or shrinking as, $\left(\frac{p}{2}+\frac{1}{n}\right)+1>n$, $\left(\frac{p}{2}+\frac{1}{n}\right)+1=n \operatorname{or}\left(\frac{p}{2}+\frac{1}{n}\right)+1<n$.  
\end{corollary}

\begin{corollary}\cite{LGDB22}\label{C:4.2.3}
    If a $n$-dimensional $\epsilon$-Kenmotsu manifold has a Codazzi type Ricci tensor and admits a conformal $\eta$-Ricci soliton $(g, \xi, \lambda, \mu)$, then 

\vskip.1in
1. if $\xi$ is spacelike then the soliton is expanding, steady or shrinking according as, $\left(\frac{p}{2}+\frac{1}{n}\right)+n>2$, $\left(\frac{p}{2}+\frac{1}{n}\right)+n=2$ or $\left(\frac{p}{2}+\frac{1}{n}\right)+n<2$; and 
\vskip.1in

2. if $\xi$ is timelike then the soliton is expanding, steady or shrinking according as, $\left(\frac{p}{2}+\frac{1}{n}\right)+2>n$, $\left(\frac{p}{2}+\frac{1}{n}\right)+2=n \operatorname{or}\left(\frac{p}{2}+\frac{1}{n}\right)+2<n$.
\end{corollary}

\begin{theorem}\cite{LGDB22}\label{T:4.2.1}
Given a $n$-dimensional $\epsilon$-Kenmotsu manifold $(M, g)$, let $(g, \xi, \lambda, \mu)$ be a conformal $\eta$-Ricci soliton. In the event that the manifold meets the curvature criterion $R(X, Y) \cdot S=0$, it can be said to be Ricci semi-symmetric. In this case, $\lambda=\left(\frac{p}{2}+\frac{1}{n}\right)+\epsilon(n-2)$ and $\mu=1$. Additionally,
\vskip.1in

1. In the event that $\xi$ is spacelike, the soliton is either expanding, steady, or shrinking as, $\left(\frac{p}{2}+\frac{1}{n}\right)>(2-n)$, $\left(\frac{p}{2}+\frac{1}{n}\right)=(2-n)$ or $\left(\frac{p}{2}+\frac{1}{n}\right)<(2-n) ;$ and
\vskip.1in

2. In the event that $\xi$ is timelike, the soliton is either expanding, steady, or shrinking as, $\left(\frac{p}{2}+\frac{1}{n}\right)+(2-n)>0$, $\left(\frac{p}{2}+\frac{1}{n}\right)+(2-n)=0$ or $\left(\frac{p}{2}+\frac{1}{n}\right)+(2-n)<0$.
\end{theorem}

\section{Solitonic inequalities on K\"ahler manifold}\label{S5}

\subsection{A finding on compactness for K\"ahler Ricci soliton}\label{S5.1}
Cao and Sesum \cite{CS07} showed results of compactness for K\"ahler Ricci gradient compact shrinking solitons. In case that $\left(M_i, g_i\right)$ represents a series of K\"ahler Ricci solitons with real-dimension  with $n \geqslant 4$, curvatures with uniformly bounded $L^{n / 2}$ norms, Ricci curvatures with uniform boundedness from below, and $\mu\left(g_i, 1 / 2\right) \geqslant A$ (where $\mu$ is the Perelman functional), A subsequence $\left(M_i, g_i\right)$ leads to a compact orbifold $\left(M_{\infty}, g_{\infty}\right)$ with finitely many isolated singularities,  where $g_{\infty}$ is a K\"ahler Ricci soliton metric in the sense of an orbifold (fulfills the K\"ahler Ricci soliton equation in a lifting about singular points and smoothly extends in some gauge to a metric that satisfies the soliton equation away from singular points.).

The authors proved the following theorem 

\begin{theorem} \cite{CS07} \label{T:5.1.1}  
Let $\left(M_i, g_i\right)$ be a series of real-dimension $n \geqslant 6$ K\"ahler Ricci solitons such that $c_1\left(M_i\right)>0$.
$$
g_i(t)-\operatorname{Ric}\left(g_i(t)\right)=\partial \bar{\partial} u_i(t)=\frac{d}{d t} g_i(t),
$$
with $\nabla_j \nabla_k u_i=\bar{\nabla}_j \bar{\nabla}_k u_i=0$, such that
(i) $ C_1 \geq \int_{M_i}|\operatorname{Rm}|^{n / 2} d V_{g_i}$,
(ii) $-C_2\leqslant\operatorname{Ric}\left(g_i\right) $ and
(iii) $A \leqslant \mu\left(g_i, 1 / 2\right)$.

For a few constants $C_1, C_2, A$ that are uniform and independent of $i$. Then, for any given orbifold $Y$ with finitely many isolated singularities, there there exists a subsequence $\left(M_i, g_i\right)$ converging to $(Y, \bar{g})$, where $\bar{g}$ is a K\"ahler Ricci soliton in the sense of an orbifold.
\end{theorem}

\subsection{On the K\"ahler-Ricci flow and Harnack's inequality}\label{S5.2}

In 1992, H.-D. Cao \cite{Cao92} studied Richard Hamilton's Ricci flow (see \cite{Ha82,Ha86}):
\begin{equation}\label{hrnk}
\frac{\partial}{\partial t} g_{i j}=-R_{i j}+g_{i j} .
\end{equation}

On a compact $n$-dimensional K\"ahler manifold $X$ with positive holomorphic bisectional curvature.The solution of \eqref{hrnk} yields the following Harnack estimate for the Ricci tensor.

\begin{theorem}\label{T:5.2}\cite{Cao92}
Given a compact K\"ahler manifold $X$ with positive bisectional curvature and $0 \leqq t<\infty$, let $g_{i j}$ be a solution of \eqref{hrnk}. Given any pair $x \in X$ and $v \in T_{x} X$, let
$$
Q_{i j}=\frac{\partial}{\partial t} R_{i j}+R_{i k} R_{k j}-R_{i j}+R_{i j, k} v^{k}+R_{i j, k} v^{\bar{k}}+R_{i j k i} v^{k} v^{i}+\frac{1}{1-e^{-t}} R_{i j}.
$$
Then, for any $t>0$, and $w \in T_{x} X, w \neq 0$, we have $Q_{i j} w^{i} w^{j}>0.$
\end{theorem}

The author obtained gradient estimates for the scalar curvature and the determinant of the Ricci tensor of the metric $g_{i j}$ by taking different traces of $Q_{i j}$.

\begin{corollary}\cite{Cao92}\label{C:5.1}
The estimate is satisfied by the scalar curvature $R$.
$$
\frac{\partial R}{\partial t}+R_{, k} v^{k}+R_{, \bar{k}} v^{\bar{k}}+R_{i j} v^{i} v^{\bar{j}}+\frac{R}{1-e^{-t}}>0
$$
for all $t>0, x \in X$ and $v \in T_{x} X$. In particular, taking $v_{i}=-R_{, i} / R$, we obtain
$$
\frac{\partial R}{\partial t}-\frac{|D R|^{2}}{R}+\frac{R}{1-e^{-t}}>0
$$
for all $t>0, x \in X$.
\end{corollary}

\begin{corollary}\cite{Cao92}\label{C:5.2}
Let $\phi=\operatorname{det}\left(R_{i j}\right) / \operatorname{det}\left(g_{i j}\right)$, then we have
$$
\frac{\partial \phi}{\partial t}-\frac{|D \phi|^{2}}{n \phi}+\frac{n \phi}{1-e^{-t}}>0
$$
for all $t>0, x \in X$.
\end{corollary}

The author got the following Harnack inequalities for the scalar curvature and the determinant of the Ricci tensor by integrating the above estimates as in \cite{LY86}.

\begin{theorem}\cite{Cao92}\label{T:5.3} Let $g_{i j}$ be the solution of \eqref{hrnk} as in \eqref{T:5.2}. Then for any $x, y \in X$ and $0<t_{1}<t_{2}<\infty$, we have:

(i) The scalar curvature $R$ satisfies the inequality
$$
R\left(x, t_{1}\right) \leqq \frac{e^{t_{2}}-1}{e^{t_{1}}-1} e^{\Delta / 4} R\left(y, t_{2}\right)
$$

(ii) The determinant of the Ricci tensor $\phi=\operatorname{det}\left(R_{i j}\right) / \operatorname{det}\left(g_{i j}\right)$ satisfies the inequality
$$
\phi\left(x, t_{1}\right) \leqq\left(\frac{e^{t_{2}}-1}{e^{t_{1}}-1}\right)^{n} e^{n \Delta / 4} \phi\left(y, t_{2}\right),
$$
where $\Delta$ is defined as
$$
\Delta=\Delta\left(x, y ; t_{1}, t_{2}\right)=\inf _{\gamma} \int_{t_{1}}^{t_{2}}\left(\left|\gamma^{\prime}(s)\right|_{s}^{2}\right) d s
$$
\end{theorem}
Taking the infimum over all curves from $x$ to $y$, where $\left|\gamma^{\prime}(s)\right|_{s}$ is the velocity of $\gamma$ at time $s$.

\section{Solitonic inequalities on compactness and non-compactness of manifold}\label{S6}

\subsection{Ricci soliton structures and non-compact manifolds as extremes of log-Sobolev inequality}\label{S6.1}

In 2019, M. Rimoldi and G. Veronelli \cite{RV19} shown that exponential decline occurs at the extremals' infinity when Ricci curvature is likewise constrained from above.
These analytical results led them to prove that a gradient Ricci soliton structure is supported by non-trivial shrinking Ricci solitons, under the same conditions.

\begin{theorem}\label{T:6.1.1}\cite{RV19}
Given a non-compact complete Riemannian manifold (connected) $\left(M^{m}, g\right)$, let us assume that
\begin{equation}
-(m-1) K \leq \operatorname{Ric} \text { and } \operatorname{inj}_{(M, g)} \geq i_{0}>0 .
\end{equation}
In the case $K \in[0,+\infty],\, i_{0} \in \mathbb{R}^{+}$, there exists a smooth extremal $v$ for the log-Sobolev functional $\mathcal{L}$ if $\lambda<\lambda_{\infty}$. Furthermore, we suppose that in place of the bound $\operatorname{Ric} \geq-(m-1) K$, assume that
$(m-1) K \geq |\operatorname{Ric}|.$
After fixing a point $o \in M$, the extremal $v$ fulfills the following positive constants: $C, c>0$ such that
\begin{equation}\label{3.10.2}
C e^{-c d^{2}(x, o)} \geq v(x) , \quad x \in M.
\end{equation}
\end{theorem}

\begin{remark}\cite{RV19}
    
 (a) The scalar curvature $R$ is contained in the log-Sobolev functional $\mathcal{L}$, however deleting the term holding the scalar curvature does not change the outcome. There isn't much to change in the proof.

(b) See \cite{Zhang12} for remarks on the condition's universality at infinity $\lambda<\lambda_{\infty}$, as well as for several examples of Riemannian manifolds that meet this condition.
\end{remark}

Rimoldi and Veronell obtained the following geometric consequence, which was also demonstrated in \cite{RV19}, by utilizing Theorem \eqref{T:6.1.1}.

\begin{theorem}\cite{RV19}
Consider a connected full non-compact Riemannian manifold $\left(M^{m}, g\right)$ that supports a Ricci soliton structure that is shrinking. Assume that $K$ and $i_{0}$ are positive constants such that
$$
(m-1) K \geq |\operatorname{Ric}| , \text { and } \operatorname{inj}_{(M, g)} \geq i_{0}>0.
$$
Moreover, $\lambda<\lambda_{\infty}$. A shrinking gradient Ricci soliton structure is then also supported by $\left(M^{m}, g\right)$.
\end{theorem}

\subsection{Theorems of compactness for gradient Ricci solitons}\label{S6.2}

X. Zhang's compactness theorem \cite{Zhang06} from 2006 requires a subsequence with uniformly bounded curvatures when a sequence of compact gradient Ricci solitons of dimension n converges to a compact orbifold with finitely many isolated singularities.

The author of \cite{Zhang06} discussed the finding of compactness for Ricci solitons in this work. When the underlying manifolds are closed (compact, boundaryless), it is simple to show that the steady and expanding Ricci solitons are really Einstein metrics. As a result, the author focused on the shrinking situation and developed the Theorem that followed.

\begin{theorem}\label{T:6.2.1}\cite{Zhang06}
 Assume that $\left(M_{\alpha}, g_{\alpha}\right)$ is a sequence of shrinking gradient Ricci solitons with dimension $n \geq 4$, meaning that it fulfills the subsequent equation:
\begin{equation}\label{CT1}
g_{\alpha}-\operatorname{Ric}\left(g_{\alpha}\right)=\nabla \mathrm{d} u_{\alpha}
\end{equation}
such that
(1) $\operatorname{Ric}\left(g_{\alpha}\right) \geq-C_{1} g_{\alpha}$,
(2) $\operatorname{diam}\left(M_{\alpha}, g_{\alpha}\right) \leq C_{2}$,
(3) $\operatorname{Vol}\left(M_{\alpha}, g_{\alpha}\right) \geq C_{3}$ and
(4) $\int_{M_{\alpha}}|\mathrm{Rm}|^{\frac{n}{2}} \mathrm{~d} V_{g_{\alpha}} \leq C_{4}$.

Given a few consistent positive constants $C_{1}, C_{2}, C_{3}, C_{4}$. Then, in the Cheeger-Gromov sense, there is a subsequence $\left(M_{\alpha}, g_{\alpha}\right)$ that converges to $\left(M_{\infty}, g_{\infty}\right)$, where $g_{\infty}$ is a Ricci soliton in the sense of an orbifold and $M_{\infty}$ is an orbifold with finitely many isolated singularities.
\end{theorem}
Moreover, if $n$ is odd, then $\left(M_{\infty}, g_{\infty}\right)$ is a smooth gradient Ricci soliton that, for sufficiently large $\alpha$, diffeomorphically maps to $M_{\alpha}$. This indicates a smooth convergence of $\left(M_{\alpha}, g_{\alpha}\right)$ (sub) to $\left(M_{\infty}, g_{\infty}\right)$.

Bounds on $\int_{M}|\mathrm{Rm}|^{2}$ are implied for Ricci solitons by lower bounds on Ricci curvature and volume, an upper bound on diameter, and a bound on $b_{2}(M)$. We possess the subsequent corollary.

\begin{corollary}\cite{Zhang06}\label{C:6.2.0}
 Given a sequence of shrinking gradient Ricci solitons of dimension 4, let $\left(M_{\alpha}, g_{\alpha}\right)$ be such that
(1) $\operatorname{Ric}\left(g_{\alpha}\right) \geq-C_{1} g_{\alpha}$,
(2) $\operatorname{diam}\left(M_{\alpha}, g_{\alpha}\right) \leq C_{2}$,
(3) $\operatorname{Vol}\left(M_{\alpha}, g_{\alpha}\right) \geq C_{3}$ and
(4) $b_{2}\left(M_{\alpha}\right) \leq C_{4}$
for a few uniformly positive constants $C_{1}, C_{2}, C_{3}, C_{4}$. Then, in the Cheeger-Gromov sense, there is a subsequence $\left(M_{\alpha}, g_{\alpha}\right)$ that converges to $\left(M_{\infty}, g_{\infty}\right)$, where $g_{\infty}$ is a Ricci soliton in the sense of an orbifold and $M_{\infty}$ is an orbifold with finitely many isolated singularities.
\end{corollary}

X. Zhang then made some deductions and derives certain estimates; specifically, $C_1$ estimates for functions $u_\alpha$ and a uniform bound for Perelman's function and the Ricci curvature. The author then calculated the Ricci soliton's $\epsilon$-regularity estimate.

\begin{proposition}\cite{Zhang06}\label{P:6.2.1}
 If $g$ is a gradient soliton that is either steady or expanding over a compact manifold $M$, then $g$ has to be an Einstein metric.

When $\lambda$ is a positive constant, from the estimate 
$\sup _{x \in M} u \leq C_8$, 
we have
\begin{equation*}
|\nabla u|^{2} \leq n \lambda+2 \lambda\left(\sup _{x \in M} u-\inf _{x \in M} u\right)-R \leq C_{10} ,
\end{equation*}
where $C_{10}$ is a constant depending only on $C_{1}, C_{2}, C_{3}$ and $\lambda$.
\end{proposition}
Using the following lemma, the author was able to determine a uniform upper bound for scalar curvature.

\begin{lemma}\cite{Zhang06}\label{C:6.2.1}
 Consider $\left(M_{\alpha}, g_{\alpha}\right).$ be a sequence of shrinking Ricci solitons $(\lambda=1)$ satisfying conditions (1), (2), and (3) of Theorem \eqref{T:6.2.1}, and $u_{\alpha}$ satisfying the constraint 
\begin{equation}\label{CT3}
(2 \pi)^{-\frac{n}{2}} \int_M \mathrm{e}^{-u} \mathrm{~d} V_g=1 .
\end{equation}
\end{lemma}
 There are positive constants $C_{4}, C_{5}$ that rely exclusively on $C_{1}, C_{2}$, and $C_{3}$, such that
\begin{equation*}
\left|u_{\alpha}\right|_{C^{1}} \leq C_{11} \;\; {\rm and}\;\;
R\left(g_{\alpha}\right) \leq C_{12} .
\end{equation*}

In the next portion, X. Zhang gave a uniform bound of Perelman's functional  (see \cite{Pe02}), $\mu\left(g, \frac{1}{2}\right)$, for a sequence of shrinking Ricci solitons ($M_{\alpha}, g_{\alpha}$) meeting conditions (1), (2) and (3) in Theorem \ref{T:6.2.1}. 

In \cite{Pe02}, Perelman proposed a functional fulfilling
\begin{equation*}
W(g, \varphi, \tau)=(4 \pi \tau)^{-\frac{n}{2}} \int_{M} \mathrm{e}^{-\varphi}\left[\tau\left(R+|\nabla \varphi|^{2}\right)+f-n\right] \mathrm{d} V_{g} 
\end{equation*}
under the constraint
\begin{equation*}\label{CT2}
(4 \pi \tau)^{-\frac{n}{2}} \int_{M} \mathrm{e}^{-\varphi} \mathrm{d} V_{g}=1 .
\end{equation*}
Then they defined the functional
$\mu(g, \tau)=\inf W(g, \cdot, \tau) ,
$
where all functions meeting the restriction \eqref{CT2} are taken over by $\inf$, and $\tau>0$.

\begin{lemma}\cite{Zhang06}\label{C:6.2.2}
 If $(M, g)$ is a shrinking gradient Ricci soliton, that is, 
$g-\operatorname{Ric}(g)=\nabla \mathrm{d} u$,
where $u$ is a minimizer of Perelman's functional $W$ with regard to metric $g$ and $\tau=\frac{1}{2}$, then $u$ satisfies the constraint \eqref{CT3}.
\end{lemma}

\begin{proposition}\cite{Zhang06}\label{P:6.2.2}
Let $u_{\alpha}$ satisfy the constraint \eqref{CT3}, and let $\left(M_{\alpha}, g_{\alpha}\right)$ be a sequence of shrinking Ricci solitons $(\lambda=1)$ satisfying conditions (1), (2), and (3) in Theorem \eqref{T:6.2.1}. Then, there is a constant $C_{6}$ that depends only on $C_{1}, C_{2}$, and $C_{3}$, such that
$\left|\mu\left(g_{\alpha}, \tfrac{1}{2}\right)\right| \leq C_{6}. $
\end{proposition}

The following compactness theorem for Ricci solitons can be obtained with ease using the Gromov-Cheeger compactness theorem.

\begin{proposition}\cite{Zhang06}\label{P:6.2.3}
 Let $\left(M_{\alpha}, g_{\alpha}\right)$ be a sequence of shrinking gradient Ricci solitons of dimension $n=3$, such that
(1) $\operatorname{Ric}\left(g_{\alpha}\right) \geq-C_{1} g_{\alpha}$,
(2) $\operatorname{diam}\left(M_{\alpha}, g_{\alpha}\right) \leq C_{2}$ and
(3) $\operatorname{Vol}\left(M_{\alpha}, g_{\alpha}\right) \geq C_{3}>0$.
Consider the following uniform constants: $C_{1}, C_{2}, C_{3}$. A smooth gradient Ricci soliton $\left(M_{\infty}, g_{\infty}\right)$ is then a subsequence in $C^{\infty}$ topology, and it converges to $\left(M_{\infty}, g_{\infty}\right)$.
\end{proposition}

The following mean value inequality was determined by the author using Moser's iteration argument.

\begin{lemma}\cite{Zhang06}\label{L:6.2.3}
     Given a compact Riemannian manifold $(M, g)$ and a Lipschitz function $f$ that satisfies
\begin{equation*}
f \triangle f \geq-\theta_{1}|\nabla f|^{2}-\theta_{2} f^{2}-\theta_{3} f^{3} 
\end{equation*}
in a feeble manner. In the event where $\theta_{1} \leq \frac{1}{4}$, a constant $\epsilon$ will exist, and it will only depend on the dimensions of $M, \theta_{3}$, and the lower bound of the Sobolev constant $C_{s}$. This means that if
$\int_{B_{P}(2 r)} f^{\frac{n}{2}} \mathrm{~d} v_{g} \leq \epsilon ,$
then
\begin{equation*}
\sup _{B_{P}\left(\frac{r}{2}\right)} f \leq C_{*}\left(1+\frac{1}{r^{2}}\right)\Big(\int_{B_{P}(r)} f^{\frac{n}{2}} \mathrm{~d} v_{g}\Big)^{\frac{2}{n}} ,
\end{equation*}
where the Sobolev constant $C_{s}$, the lower bound of $\operatorname{Vol}(M)$, and the dimensions of $M, \theta_{2}, \theta_{3}$ are the only factors that affect $C_{*}$.
\end{lemma}

\subsection{Geometry of compact shrinking Ricci solitons}\label{S6.3}
In 2014, B.-Y. Chen and S. Deshmukh \cite{chendesh} proved two characterizations of compact shrinking trivial Ricci solitons.
They proved the following results
\begin{theorem}\label{T:3} Let $(M,g,f,\lambda )$ be an $n$-dimensional compact connected shrinking gradient Ricci soliton of positive Ricci curvature. If the Ricci curvature $Ric$ and the scalar curvature $S$ of $(M,g)$ satisfy%
\begin{equation*}
Ric(\nabla S,\nabla S)\leq \lambda \left( \left\Vert \nabla S\right\Vert
^{2}+\frac{\lambda }{2}\left( n^{2}\lambda ^{2}-S^{2}\right) \right) \text{,}
\end{equation*}%
then $S$ is a solution of the Poisson equation
$\Delta \varphi =\sigma $ with $\sigma =\lambda (n\lambda -S)$.\end{theorem}

\begin{corollary}\label{C:2} A compact connected shrinking gradient Ricci soliton $(M,g,f,\lambda )$ of positive Ricci curvature with $\lambda<\lambda_{1}$ is trivial if and only if the scalar curvature $S$ satisfies%
\begin{equation*}
Ric(\nabla S,\nabla S)\leq \lambda \left( \left\Vert \nabla S\right\Vert^{2}+\frac{\lambda }{2}\left( n^{2}\lambda ^{2}-S^{2}\right) \right),\quad n=\dim M \text{.} 
\end{equation*}
\end{corollary}

On a compact Riemannian manifold $(M,g)$ and a function $\varphi :M\rightarrow R$, the average value of $\varphi $ is a real number defined by%
\begin{equation*} \varphi_{av}=\frac{1}{{\rm Vol}(M)}\int_{M}\varphi .\end{equation*}%

Chen and Deshmukh also proved the following characterization of trivial Ricci solitons.

\begin{theorem}\label{T:1} \cite{chendesh}  A compact connected shrinking gradient Ricci soliton $(M,g,f,\lambda )$ with normalized potential is trivial if and only if
 $ (fS)_{av}\leq \frac{1}{2}n^2\lambda$,
where $n=\dim M$ and $S$ is the  scalar curvature of $(M,g)$.\end{theorem}

\section{Solitonic inequalities on curvature}\label{S7}

\subsection{About the Ricci curvature of Lagrangian and isotropic submanifolds in complex space forms}\label{S7.1}
B.-Y. Chen proved in \cite{Chen00} that any isotropic submanifold $M^{n}$ in a complex space form $\tilde{M}^{m}(4 c)$ satisfies
\begin{equation}\label{3.51}
 (n-1) c+\frac{n^{2}}{4} H^{2} \geqq \overline{R i c}  .
\end{equation}
In \cite{Chen00}  he also proved that a Lagrangian submanifold of a complex space form satisfying the equality case of \eqref{3.51} is a minimal submanifold. Finally, he gave an explanation of the geometry of Lagrangian submanifolds satisfying \eqref{3.51} for equality, as long as the second basic form kernel's dimension doesn't change.
\begin{theorem} \cite{Chen00}
     Given a complex space form $\tilde{M}^{m}(4 c)$ and an isotropic submanifold $M^{n}$, the Ricci tensor $S$ of $M^{n}$ satisfies
\begin{equation}\label{3.52}
\left(\frac{n^{2}}{4} H^{2} +(n-1) c\right)g \geqq S.
\end{equation}
The equality sign of \eqref{3.52} holds ientically if and only if $M^{n}$ is a fully geodesic submanifold or if $n=2$ and $M^{n}$ is entirely umbilical.\\
\end{theorem}

Then he proved the main result on minimality of Lagrangian in the following Theorem
\begin{theorem}\label{LagR2} \cite{Chen00}
Let $M^{n}$ be a Lagrangian submanifold of $\tilde{M}^{n}(4 c)$, which is a complex space form. Then
\begin{equation}\label{3.53}
\overline{R i c} \leqq(n-1) c+\frac{n^{2}}{4} H^{2} .
\end{equation}
\end{theorem}

Clearly, Theorem \eqref{LagR2} follows  from the following.

\begin{lemma}\cite{Chen00}
We have \eqref{3.53} for every n-dimensional isotropic submanifold of a complex space form $\tilde{M}^{m}(4 c)$.
     Furthermore, if an isotropic submanifold $M^{n}$ of $\tilde{M}^{m}(4 c)$ satisfies the equality case of \eqref{3.53} at a point $p$, then $J\left(T_{p} M^{n}\right)$ is perpendicular to the mean curvature vector $\vec{H}$ at $p$.
\end{lemma}

\subsection{Ricci solitons with finite scalar curvature ratio and complete gradient expansion asymptotically}\label{S7.2}

Given $\left(M^{n}, g, f\right), n \geq 5$ is a complete expanding gradient  Ricci soliton with nonnegative Ricci curvature.
In 2023, Cao et al., \cite{CLX23}, showed that if the scalar asymptotic curvature ratio of $\left(M^{n}, g, f\right)$ is finite (i.e., $\left.\operatorname{lim~sup}_{r \rightarrow \infty} R r^{2}<\infty\right)$, then the tensor of  Riemann curvature must have at least sub-quadratic decay, i.e., $\lim _{\sup _{r \rightarrow \infty}}|R m| r^{\alpha}<\infty$ for any $0<\alpha<2$.
The following are the main result given by them 
\begin{theorem}\cite{CLX23} 
 A $n$-dimensional expanding complete gradient  Ricci soliton with finite asymptotic scalar curvature ratio and nonnegative Ricci curvature $R c \geq 0$ is given. Let $\left(M^{n}, g, f\right), n \geq 5$. We have
\begin{equation*}\label{3.91}
\limsup _{r \rightarrow \infty} R r^{2}<\infty ,
\end{equation*}
where the distance function to a given base point $x_{0} \in M$ is denoted by $r=r(x)$. 
For any $0<\alpha<2$, $\left(M^{n}, g, f\right)$ has a finite $\alpha$-asymptotic curvature ratio:
\begin{equation*}\label{3.92}
A_{\alpha}:=\limsup _{r \rightarrow \infty}|R m| r^{\alpha}<\infty . 
\end{equation*}
Additionally, given $n$ and the geometry of $\left(M^{n}, g, f\right)$, there are constants $C>0$ that determine the sequences $\left\{r_{j}\right\} \rightarrow \infty$ and $\left\{\alpha_{j}\right\} \rightarrow 2$ such that
$$
|R m|(x) \leq C(r(x)+1)^{-\alpha_{j}}.
$$
For each $x \in M \backslash B\left(x_{0}, r_{j}+1\right)$.
Where B represents a geodesic ball.$f\in C^\infty(M)$, and A represents asymptotic curvature.
\end{theorem}

\subsection{Estimates of curvature for steady Ricci solitons}\label{S7.3}

P.-Y. Chan's work \cite{Chan19} from 2019 improves the estimate given in \cite{CC20} for an $n$-dimensional full non-Ricci flat gradient stable Ricci soliton with a potential function constrained by a constant and a curvature tensor $R m$ fulfilling $R m|<5$, $|R m| < C e^{-r}$ for some constant $C>0$. For a four-dimensional full non-Ricci flat gradient steady Ricci soliton, the decay rate of $|R m|$ is exponential.\\

In \cite{MSW19}, O. Munteanu, C.-J.A.  Sung and J. Wang investigated the solvability of the weighted Poisson equation for a certain class of smooth metric measure spaces. They illustrated the following, for example:\\

\begin{theorem}
 \cite{MSW19} Given $n \geq 2$, let $\left(M^{n}, g, f\right)$ be a $n$ dimensional complete non Ricci flat gradient steady Ricci soliton. In the event that a constant limits the potential function $f$ above and $\lim _{r \rightarrow \infty} r|R m|=0$, then a positive constant $C$ exists such that
\begin{equation}\label{CE1}
|R m|(x) \leq C(1+r(x))^{3(n+1)} e^{-r(x)} 
\end{equation}
\end{theorem}
where $r=r(x)$ is the distance of $x$ from a fixed point $p_{0} \in M$.\\

Instead of using Green's function estimate in \cite{MSW19}, the author of \cite{Chan19} sharpens the upper bound under weaker curvature decay conditions by studying the curvature features of growing gradient Ricci solitons using Deruelle's maximum principle \cite{Der17}.

\begin{theorem}\label{CES1}\cite{Chan19}
 Let $\left(M^{n}, g, f\right)$ be a $n$ dimensional non-Ricci flat complete gradient steady Ricci soliton, where $n \geq 2$. Let us assume that the curvature tensor $R m$ satisfies $\underset{r \rightarrow \infty}{\lim \sup } r|R m|<\frac{1}{5}$ and that the potential function $f$ is bounded above by a constant. Next, a positive constant C exists such that
\begin{equation}\label{CE2}
|R m|(x) \leq C e^{-r(x)} \text { on } M, 
\end{equation}
where $r=r(x)$ is the distance of $x$ from a fixed point $p_{0} \in M$.
\end{theorem}

\begin{remark}\cite{Chan19}
When $\Sigma$ and $\mathbb{T}^{n-2}$ represent the Hamilton's cigar soliton and the $n-2$ dimensional flat torus, respectively, the decay rate is sharp on $\Sigma \times \mathbb{T}^{n-2}$.\\
\end{remark} 

\begin{remark}\cite{Chan19}
In Theorem \eqref{CES1}, there exists a technical constant $\frac{1}{5}$. According to \cite{Bre13} and \cite{Bre14}, $\lim _{r \rightarrow \infty} r|R m|=\sqrt{\frac{(n-1)}{2(n-2)}}>\frac{1}{5}$ for the $n$ dimensional Bryant soliton. Whether the constant $\frac{1}{5}$ is optimal is not obvious at this time.
\end{remark}

\begin{remark}\cite{Chan19}   
For simplicity's sake, the author states the Theorem as follows: $\varlimsup_{r \rightarrow \infty} r^{-1} f<1$. This condition on $f$ seems to be weaker than the previous one, but it does not affect the validity of the conclusion \eqref{CE2}. For $f$ to be bounded above by a constant, it is necessary that ric $\geq 0$ and $S \rightarrow 0$ at infinity (see \cite{CN09}).\\
\end{remark}

O. Munteanu and J. Wang \cite{MW15} noted, among other unique characteristics, that the Riemann curvature $R m$ in a four-dimensional gradient Ricci soliton can be bounded by $\nabla$\! Ric and Ric.
\begin{equation}\label{3.163}
|R m| \leq A_{0}\left(|\operatorname{Ric}|+\frac{|\nabla \mathrm{Ric}|}{|\nabla f|}\right) 
\end{equation}
for some universal positive constant $A_{0}$.\\

Here is another main result of this paper:
\begin{theorem}\cite{Chan19}
Let $\left(M^{4}, g, f\right)$ be a four-dimensional non-Ricci flat gradient steady Ricci soliton with $\lim _{r \rightarrow \infty} S=0$. There is a positive
constant $c$ such that
\begin{equation}
|R m| \leq c S \; \text { on } M.
\end{equation}
\end{theorem}

A 4-dimensional full nontrivial gradient stable Ricci soliton's Riemann curvature $R m$ decays exponentially if the potential $f$ is constrained from above and $\overline{\lim }_{r \rightarrow \infty} r S$ is small enough, as was shown in \cite{Chan19}.
\begin{theorem}\cite{Chan19}
Let $\left(M^{4}, g, f\right)$ be a four-dimensional non-Ricci flat gradient steady Ricci soliton with $\lim _{r \rightarrow \infty} S=0$. Assume the potential function $f$ is bounded from above by a constant and $\varlimsup_{r \rightarrow \infty} r S<\frac{1}{5 A_{0}^{2}}$, where $A_{0}$ is the constant in \eqref{3.163}. Then there exists a constant $C>0$ such that
$|R m|(x) \leq C e^{-r(x)}$ on $M.$
\end{theorem}

\section{Solitonic inequalities on stability and instability}\label{S8}

\subsection{ stability and instability of Ricci soliton} 

  Kr\"oncke \cite{Kr15} investigated the volume-normalized Ricci flow near compact shrinking Ricci Solitons. Any normalized Ricci flow that begins near a compact Ricci Soliton $(M, g)$, which is a local maximum of Perelman's shrinker entropy, is shown to exist for all time and to converge towards a Ricci Soliton. If $g$ is not a local maximum of the shrinker entropy, then a nontrivial normalized Ricci flow emerges from it.
\begin{lemma}\cite{Kr15}
    Consider a gradient shrinking Ricci soliton $\left(M, g_{0}\right)$. Then, in the space of metrics, there is a constant $C>0$ and a $C^{2, \alpha_{-}}$ neighbourhood $\mathcal{U}$ of $g_{0}$ such that
\begin{align*}
\left\|\left.\frac{d}{d t}\right|_{t=0} f_{g+t h}\right\|_{C^{2, \alpha}} & \leq C\|h\|_{C^{2, \alpha}}, \quad\left\|\left.\frac{d}{d t}\right|_{t=0} f_{g+t h}\right\|_{H^{i}} \leq C\|h\|_{H^{i}},  \;\; 
\left.\left|\frac{d}{d t}\right|_{t=0} \tau_{g+t h} \right\rvert\,  \leq C\|h\|_{L^{2}}, \;\;i=0,1,2, 
\end{align*}
for  all $g \in \mathcal{U}$ ,
where g, $g_{0}$ are Riemannian metric tensors, $f \in C^{\infty}\left(M\right)$, and H, L are the maps defined in \cite{Kr15}.
\end{lemma}

\begin{theorem} \cite{Kr15}\label{T:8.1}
    Consider a Ricci soliton that shrinks gradient shrinkage, $\left(M, g_{0}\right)$. Then, given $g_{0}$ and constants $\sigma \in[1 / 2,1), C>0$, there exists a $C^{2, \alpha}$ neighbourhood $\mathcal{U}$ such that
    \begin{equation}
C\left\|\tau\left(\operatorname{Ric}_{g}+\nabla^{2} f_{g}\right)-\frac{1}{2} g\right\|_{L^{2}}\geq\left|v(g)-v\left(g_{0}\right)\right|^{\sigma} 
\end{equation}
for all $g \in \mathcal{U}$.
\end{theorem}

Further, they proved the main results on dynamical stability and instability by considering the $\tau$-flow 
\begin{equation}\label{3.33}
\dot{g}(t)=-2 \operatorname{Ric}_{g(t)}+\frac{1}{\tau_{g(t)}} g(t) 
\end{equation}
It is well defined around a gradient-shrinking Ricci soliton. Note that under the $\tau$-flow, $v$ is nondecreasing.The writers built a modified $\tau$-flow in this way: The collection of diffeomorphisms produced by $X(t)=-\operatorname{grad}_{g(t)} f_{g(t)}$ is denoted as $\tilde{g}(t):=\varphi_{t}^{*} g(t)$, where $g(t)$ is a solution of \eqref{3.33}. Let $\varphi_{t} \in \operatorname{Diff}(M), t \geq 1$ be the collection of diffeomorphisms. Then 
\begin{equation}
\frac{d}{d t} \tilde{g}(t)=-2\left(\operatorname{Ric}_{\tilde{g}(t)}+\nabla^{2} f_{\tilde{g}(t)}\right)+\frac{1}{\tau_{\tilde{g}(t)}} \tilde{g}(t). 
\end{equation}
This represents the $\tau$ gradient flow in relation to the weighted $L^{2}$-measure.
\begin{lemma} \cite{Kr15}
($\tau$-flow Shi estimates) Assuming that $g(t), t \in[0, T]$ is a solution to the $\tau$-flow \eqref{3.33},
$$
T^{-1} \geq \quad \sup _{p \in M}\left|R_{g(t)}\right|_{g(t)}+\frac{1}{\tau_{g(t)}} \quad \forall t \in[0, T].$$
Then, there exists a constant $C(k)$ such that, for each $k \geq 1$, we have
$$\sup _{p \in M}\left|\nabla^{k} R_{g(t)}\right|_{g(t)} \leq C(k) \cdot T^{-1} t^{-k / 2} \quad \forall t \in(0, T].$$
\end{lemma}

\begin{theorem} \cite{Kr15}\label{T:8,2}
    (Dynamical stability) Assume that $k \geq 3$ and that the gradient decreasing Ricci soliton is $(M, g)$. Suppose that a local maximizer of $v$ is $g$. Then for every $C^{k}$-neighbourhood $\mathcal{U}$ of $g$, there is a $C^{k+5}$-neighbourhood $\mathcal{V}$ for which the following is true:
    A 1-parameter family of diffeomorphisms $\varphi_{t}$ exists for every metric $g_{0} \in \mathcal{V}$, such that the modified flow $\varphi_{t}^{*} g(t)$ always stays in $\mathcal{U}$ and converges to a gradient shrinking Ricci soliton $g_{\infty}$ in $\mathcal{U}$ as $t \rightarrow \infty$ for the $\tau$-flow \eqref{3.33} beginning at $g_{0}$. Polynomial rate convergence suggests that there are constants $C, \alpha>0$, such that
$$
\left\|\varphi_{t}^{*} g(t)-g_{\infty}\right\|_{C^{k}} \leq C(t+1)^{-\alpha} .
$$
\end{theorem}

\begin{lemma} \cite{Kr15}
   Consider a gradient-shrinking Ricci soliton $\left(M, g_{0}\right)$. Then, for any $g_{0}$, there is a $C^{2, \alpha}$ neighbourhood $\mathcal{U}$ with a constant $C>0$, such that
$$
\begin{aligned}
  C\|k\|_{C^{2, \alpha}}\|h\|_{H^{1}} & \geq\left|\frac{d^{2}}{d t d s}\right|_{t, s=0} f_{g+s k+t h} \|_{H^{1}}\;\; {\rm and}\;\;
C\|k\|_{C^{2, \alpha}}\|h\|_{H^{1}} & \geq \left.\left|\frac{d^{2}}{d t d s}\right|_{t, s=0} \tau_{g+s k+t h} \right\rvert\, , 
\end{aligned}
$$
for all $g \in \mathcal{U}$.
\end{lemma}
\begin{proposition}\cite{Kr15}
     (Estimates of the third variation of $v)$ Consider a gradient shrinking Ricci soliton $\left(M, g_{0}\right)$. It is possible for $g_{0}$ to have a $C^{2, \alpha}$-neighborhood $\mathcal{U}$ such that
$$
\left.\left|\frac{d^{3}}{d t^{3}}\right|_{t=0} v(g+t h) \right\rvert\, \leq C\|h\|_{H^{1}}^{2}\|h\|_{C^{2, \alpha}}
$$
for all $g \in \mathcal{U}$ and some constant $C>0$.
\end{proposition}

\section{Solitonic inequalities on warped product manifolds}\label{S9}

For general references on warped product manifolds we refer to the books \cite{book2017,Obook83}.

\subsection{Inequalities of Ricci soliton in CR-warped product manifold}\label{S9.1}

In 2023, Y. Li et al., \cite{LSMKA23}, applied inequalities by using the notion of CR- warped products introduced in \cite{Chen01}. 

\begin{theorem} \cite{LSMKA23}
     Assume that $M=M_{T} \times_{f} M_{\perp}$ is a $\mathcal{C} \mathcal{R}$-warped product in $\tilde{M}(c)$ that admits a Ricci soliton with a shrinking gradient. Thus, the following disparity is true:
\begin{equation}
\|h\|^{2} \geq \beta\|\nabla(\ln f)\|^{2}+c \alpha \beta+\beta \sum_{q=1}^{2 \alpha} \operatorname{Ric}\left(E_{q}, E_{q}\right).
\end{equation}
Furthermore, the equality holds if and only if $M_{\perp}$ is a fully umbilical submanifold of $M(c)$ and $M_{T}$ is a totally geodesic submanifold.
\end{theorem}

\subsection{Ricci solitons on noncompact warped products}\label{S9.2}

Ricci solitons with complete noncompact warped product gradient were studied by V. Borges and K. Tenenblat \cite{BT22}. It is proved hat the warping function's gradient and nonexistence results are true. When the soliton is expanding or steady, some pde estimates and rigidity discovered when investigating warped product Einstein manifolds are applied to a wider context. A nonexistence theorem is presented for shrinking soliton; this theorem has no counterpart in the Einstein case and is demonstrated by utilizing the properties of a weighted Laplacian's first eigenvalue.The noncompact warped product Ricci soliton yielded the following results.

\begin{theorem}\cite{BT22}
     Let $M^n \times_h F^m$ be a complete gradient shrinking Ricci soliton and $R_F$ the scalar curvature of $F^m$. Suppose that $M$ is noncompact, fix $q_0 \in F$ and let $\mu$ be defined by $m \mu=R_F\left(q_0\right)$. The following are true:

\vskip.06in     
1. If $h$ is not constant, then $\mu>0$.

\vskip.06in     
2. If $h \leq \sqrt{\frac{\mu}{\lambda}}$, then $h$ is constant.
\end{theorem}

\begin{theorem}\cite{BT22}
Consider a complete gradient steady Ricci soliton $M^n \times_h F^m$, where the scalar curvature of $F^m$ is denoted by $R_F$. Assume that $M$ is noncompact. Fix $q_0 \in F$. Let $m \mu=$ $R_F\left(q_0\right)$ define $\mu$. The following  are true:

\vskip.06in     
 1. If $h$ is not constant, then $\mu>0$.

 \vskip.06in     
2. If $\sup h<+\infty$, then $h$ is constant.
 \end{theorem}
 
\begin{theorem}\cite{BT22}
    Assume that $M^n \times_h F^m$ is a complete gradient expanding Ricci soliton and that $R_F$ is the scalar curvature of $F^m$. Let's say $M$ is not compact. Adjust $q_0 \in F$. To define $\mu$, let $m \mu=R_F\left(q_0\right)$. The following are true:
    
\vskip.06in     
1. If $\mu<0$, then $h \geq \sqrt{\frac{\mu}{\lambda}}$,

\vskip.06in     
2. If $\mu \leq 0$, then $|\nabla \ln h|^2 \leq-\frac{\lambda}{m}$,

\vskip.06in     
3. If $\mu<0$ and $\sup _M h<+\infty$, then $|\nabla \ln h|^2 \leq-\frac{\lambda}{m}+\frac{2 \mu}{m\left(\sup _M h\right)^2}$,

\vskip.06in     
4. If $\mu \geq 0$ and $\sup h<+\infty$, then $h$ is constant.
\end{theorem}

\section{Solitonic inequalities on statistical submersion}\label{S10}

\subsection{Solitons and sharp inequalities for statistical submersion}\label{S10.1}
In 2023 Siddiqui et al. studied in  \cite{CSS24}  Ricci and scalar curvatures for a given statistical submersion and constructed some inequalities involving the Ricci and scalar curvatures:
\begin{align}\label{3.41}
\operatorname{Ric}(E, F)= & \overline{\operatorname{Ric}}(E, F)-g\left(\mathcal{T}_{E} F, N^{*}\right)+(\hat{\delta} \mathcal{T})(E, F)  +g\left(\mathcal{A} E, \mathcal{A}^{*} F\right)-g\left(\nabla_{E}^{*} \sigma, F\right) 
\end{align}
and
\begin{align}\label{3.42}
\operatorname{Ric}(X, Y)= & \hat{\operatorname{Ric}}(X, Y)+g\left(\nabla_{X}^{*} N^{*}, Y\right)-g(\mathcal{T} X, \mathcal{T} Y)+(\hat{\delta} \mathcal{A})(X, Y)\nonumber \\
& +g\left(\sigma, \mathcal{A}_{X} Y\right)-g\left(\mathcal{A}_{X}, \mathcal{A}_{Y}^{*}\right)-g\left(\mathcal{A}_{X}^{*}, \mathcal{A}_{Y}^{*}\right) 
\end{align}
where
$$
\begin{gathered}
(\hat{\delta} \mathcal{T})(E, F)=\sum_{i=1}^{n} g\left(\left(\nabla_{X_{i}} \mathcal{T}\right)(E, F), X_{i}\right), \;\; 
(\hat{\delta} \mathcal{A})(X, Y)=\sum_{j=1}^{m} g\left(\left(\nabla_{E_{j}} \mathcal{A}\right)(X, Y), E_{j}\right) \\
\sum_{i=1}^{n} g\left(\mathcal{A}_{X} X_{i}, \mathcal{A}_{Y} X_{i}\right)=g\left(\mathcal{A}_{X}, \mathcal{A}_{Y}\right)=\sum_{j=1}^{m} g\left(\mathcal{A}_{X}^{*} E_{j}, \mathcal{A}_{Y}^{*} E_{j}\right), \sigma=\sum_{i=1}^{n} \mathcal{A}_{X_{i}} X_{i}, \\
\sum_{i=1}^{n} g\left(\mathcal{A}_{X_{i}} E, \mathcal{A}_{X_{i}} F\right), \quad g(\mathcal{T} X, \mathcal{T} Y)=g(\mathcal{A} E, \mathcal{A} F)=\sum_{j=1}^{m} g\left(\mathcal{T}_{E_{j}} X, \mathcal{T}_{E_{j}} Y\right) .
\end{gathered}
$$
in which $\{E_{j}\}$ is an orthonormal basis, $\mathcal{T}$ and $\mathcal{A}$ are O'Neill tensors, and $\delta$(X) is horiontal divergence.\\

Siddiqui et al.   proved in \cite{CSS24} the following results using \eqref{3.41} and \eqref{3.42}.

\begin{theorem}\cite{CSS24}\label{T:10.1}
Let $\psi:(M, \nabla, g) \longrightarrow(N, \hat{\nabla}, \hat{g})$ be a statistical submersion. Then, we have
\begin{equation}\label{3.43}
\operatorname{Ric}(E, E) \geq \overline{\operatorname{Ri}}(E, E)-m^{2} g\left(\mathcal{T}_{E} E, H^{*}\right)+(\hat{\delta} \mathcal{T})(E, E)-(\bar{\delta} \sigma)(E, E) .
\end{equation}
\end{theorem} 
The equality case holds in \eqref{3.43} if and only if $\mathcal{H}(M)$ is integrable.\\

Since $2 \mathcal{A}^{0}=\mathcal{A}+\mathcal{A}^{*}$, we have
\begin{theorem}\cite{CSS24}\label{T:10.2}
Let $\psi:(M, \nabla, g) \longrightarrow(N, \hat{\nabla}, \hat{g})$ be a statistical submersion. Then, we have
\begin{align}\label{3.44}
\operatorname{Ric}(X, X) \leq & g\left(\nabla_{X}^{*}N^{*}, X\right)+\hat{\operatorname{Ric}}(X,X)+(\hat{\delta}\mathcal{A})(X, X) \nonumber \\
& +g\left(\sigma,\mathcal{A}_{X} X\right)-2 g\left(\mathcal{A}_{X}^{0}, \mathcal{A}_{X}^{*}\right) .
\end{align}
\end{theorem} 
The equality case holds in \eqref{3.44} if and only if each fiber is totally geodesic with respect to $\nabla(\mathcal{T}=0)$.

\begin{theorem}\cite{CSS24}\label{T:10.3}
    Let $\psi:(M, \nabla, g) \longrightarrow(N, \hat{\nabla}, \hat{g})$ be a statistical submersion. Then, we have
\begin{equation}\label{3.45}
2 R \geq 2 \bar{R}-m^{2} g\left(H, H^{*}\right)  .
\end{equation}
\end{theorem}
In \eqref{3.45}, the equality case is satisfied if and only if $\mathcal{T}$ or $\mathcal{T}^{*}$ is a multiple of the other. Specifically, every fiber is either completely geodesic with regard to $\nabla(\mathcal{T}=0)$ or completely geodesic with regard to $\nabla^{*}\left(\mathcal{T}^{*}=0\right)$.

\begin{theorem}\cite{CSS24} \label{T:10.4}   
 Let $\psi:(M, \nabla, g) \longrightarrow(N, \hat{\nabla}, \hat{g})$ be a statistical submersion. Then, we have
\begin{equation}\label{3.46} 2 R \leq 2 \hat{R}+g(\sigma, \sigma) .\end{equation}
The equality case holds in \eqref{3.46} if and only if $\mathcal{A}_{\mathcal{H}} \mathcal{H}=0$. where H is map defined in \cite{CSS24}.
\end{theorem}

Taking into account relations \eqref{3.41} and \eqref{3.42}, they derived the following equation .
\begin{align}\label{3.47}
R-\bar{R}-\hat{R}= & -2 g(\mathcal{A}, \mathcal{A})+g\left(\mathcal{A}, \mathcal{A}^{*}\right)-g\left(\mathcal{T}, \mathcal{T}^{*}\right)-g\left(N, N^{*}\right)  -\hat{\delta} N-\hat{\delta}^{*} N^{*}-\bar{\delta} \sigma+\bar{\delta}^{*} \sigma+g(\sigma, \sigma) 
\end{align}
The scalar curvatures of the vertical and horizontal spaces of $M$ are denoted by $\bar{R}$ and $\hat{R}$. Here
$$
\begin{gathered}
g\left(\mathcal{T}, \mathcal{T}^{*}\right)=\sum_{i=1}^{n} g\left(\mathcal{T} X_{i}, \mathcal{T}^{*} X_{i}\right), \;\; g(\mathcal{A}, \mathcal{A})=\sum_{i=1}^{n} g\left(\mathcal{A}_{X_{i}}, \mathcal{A}_{X_{i}}\right), \;\;
g\left(\mathcal{A}, \mathcal{A}^{*}\right)=\sum_{i=1}^{n} g\left(\mathcal{A}_{X_{i}}, \mathcal{A}_{X_{i}}^{*}\right) .
\end{gathered}
$$

By using Cauchy-Buniakowski-Schwarz inequality and equation \eqref{3.47}, we have the following Theorem.

\begin{theorem}\cite{CSS24}\label{T:10.5}
    Let $\psi:(M, \nabla, g) \longrightarrow(N, \hat{\nabla}, \hat{g})$ be a statistical submersion. Then, we have
\begin{align}\label{3.48}
R \geq & \bar{R}+\hat{R}-2\|\mathcal{A}\|^{2}+g\left(\mathcal{A}, \mathcal{A}^{*}\right)-\|\mathcal{T}\|\left\|\mathcal{T}^{*}\right\|-g\left(N, N^{*}\right) 
 -\hat{\delta} N-\hat{\delta}^{*} N^{*}-\bar{\delta} \sigma+\bar{\delta}^{*} \sigma+g(\sigma, \sigma).
\end{align}
\end{theorem}

Note that in \eqref{3.48}, the equality case is satisfied if and only if $\mathcal{T}$ or $\mathcal{T}^{*}$ is a multiple of the other. Specifically, every fiber is either completely geodesic with regard to $\nabla(\mathcal{T}=0)$ or completely geodesic with regard to $\nabla^{*}\left(\mathcal{T}^{*}=0\right)$.

\vskip.05in

\noindent{\small B.-Y. CHEN\\
ADDRESS: Department of Mathematics, Michigan State University, 619 Red Cedar Road, East Lansing, Michigan 48824--1027, U.S.A.\\
E-MAIL: chenb@msu.edu\\
{\bf ORCID ID: 0000-0002-1270-094X}\\
\\
M. A. CHODHARY \\
ADDRESS: Department of Mathematics, Maulana Azad National Urdu University, Hyderabad,
India\\
E-MAIL: majid\_alichoudhary@yahoo.co.in\\
{\bf ORCID ID: 0000-0001-5920-1227}\\
\\
M. NISAR\\
ADDRESS: Department of Mathematics, Maulana Azad National Urdu University, Hyderabad, 
India\\
E-MAIL: nisar.msc@gmail.com\\
{\bf ORCID ID: 0009-0006-0901-6123}\\
\\M. D. SIDDIQI\\
ADDRESS: Department of Mathematics, Jazan University, P.O. Box 114, Jazan 45142, Saudi
Arabia.\\
E-MAIL: msiddiqi@jazanu.edu.sa}\\
{\small {\bf ORCID ID: 0000-0002-1713-6831}}

\end{document}